\documentclass{article}
\usepackage{latexsym,amsfonts,amsmath,amsthm,makeidx,extarrows,chemarrow}
\usepackage{makeidx}
\makeindex
\newtheorem{definition}{Definition}[section]
\newtheorem{theorem}{Theorem}[section]
\newtheorem{lemma}{Lemma}[section]
\newtheorem{corollary}{Corollary}[section]

\newtheorem{remark}{Remark}[section]
\newcommand{\RN}{\mathbb R^N}
\newcommand{\Om}{\Omega}

\newcommand{\iy}{\infty}

\newcommand{\s}{\section}
\newcommand{\dd}{\delta}

\newcommand{\la}{\lambda}

\newcommand{\R}{\mathbb R}
\newcommand{\al}{\alpha}

\newcommand{\wH}{\widetilde H}

\newcommand{\bb}{\beta}

\newcommand{\e}{\varepsilon}
\newcommand{\vp}{\varphi}

\newcommand{\bt}{\begin{theorem}}
\newcommand{\et}{\end{theorem}}
\newcommand{\bl}{\begin{lemma}}
\newcommand{\el}{\end{lemma}}
\newcommand{\bd}{\begin{definition}}
\newcommand{\ed}{\end{definition}}
\newcommand{\bc}{\begin{corollary}}
\newcommand{\ec}{\end{corollary}}
\newcommand{\bp}{\begin{proof}}
\newcommand{\ep}{\end{proof}}
\newcommand{\bx}{\begin{example}}
\newcommand{\ex}{\end{example}}
\newcommand{\bi}{\begin{exercise}}
\newcommand{\ei}{\end{exercise}}
\newcommand{\bo}{\begin{prop}}
\newcommand{\eo}{\end{prop}}
\newcommand{\br}{\begin{remark}}
\newcommand{\er}{\end{remark}}
\newcommand{\be}{\begin{equation}}
\newcommand{\ee}{\end{equation}}
\newcommand{\ba}{\begin{align}}
\newcommand{\ea}{\end{align}}
\newcommand{\bn}{\begin{enumerate}}
\newcommand{\en}{\end{enumerate}}
\newcommand{\bg}{\begin{align*}}
\newcommand{\bcs}{\begin{cases}}
\newcommand{\ecs}{\end{cases}}

\newcommand{\ga}{\gamma}
\newcommand{\sg}{\sigma}
\newcommand{\bean}{\begin{eqnarray*}}
\newcommand{\eean}{\end{eqnarray*}}

\numberwithin{equation}{section}

\begin{document}

\title{\bf Multiple sign-changing and semi-nodal solutions for coupled
Schr\"{o}dinger equations\thanks{Chen and Zou are supported by NSFC (11025106).  E-mail:
chenzhijie1987@sina.com (Chen);\quad cslin@math.ntu.edu.tw (Lin); \quad
\quad wzou@math.tsinghua.edu.cn (Zou)}}
\date{}
\author{{\bf Zhijie Chen$^{1}$,   Chang-Shou Lin$^{2}$,
  Wenming Zou$^{3}$}\\
\footnotesize {\it $^{1, 3}$Department of Mathematical Sciences, Tsinghua
University,}\\
\footnotesize {\it Beijing 100084, China}\\
\footnotesize {\it $^2$Taida Institute for Mathematical Sciences, Center for Advanced Study in Theoretical Science,}\\
\footnotesize {\it National Taiwan
University, No.1, Sec. 4, Roosevelt Road, Taipei 106, Taiwan} }

\maketitle \vskip0.10in
\begin{center}
\begin{minipage}{120mm}
\begin{center}{\bf Abstract}\end{center}

We study the following coupled
Schr\"{o}dinger equations which have appeared as several models from mathematical physics:
\begin{displaymath}
\begin{cases}-\Delta u_1 +\la_1 u_1 =
\mu_1 u_1^3+\beta u_1 u_2^2, \quad x\in \Omega,\\
-\Delta u_2 +\la_2 u_2 =\mu_2 u_2^3+\beta u_1^2 u_2,  \quad   x\in
\Om,\\
u_1=u_2=0  \,\,\,\hbox{on \,$\partial\Om$}.\end{cases}\end{displaymath}
Here $\Om\subset\RN (N=2, 3)$ is a smooth bounded domain,
$\la_1, \la_2$, $\mu_1, \mu_2$ are all positive constants.
We show that, for each $k\in\mathbb{N}$ there exists $\bb_k>0$ such that
this system has at least $k$ sign-changing solutions (i.e., both two components change sign) and $k$ semi-nodal
solutions (i.e., one component changes
sign and the other one is positive) for each fixed $\bb\in (0, \bb_k)$.
\end{minipage}
\end{center}

\vskip0.1in \s{Introduction}

In this paper we study solitary wave solutions of the coupled Gross-Pitaevskii equations (cf. \cite{CLLL}):
\be\label{eq1}
\begin{cases}-i \frac{\partial}{\partial t}\Phi_1=\Delta \Phi_1+
\mu_1 |\Phi_1|^2 \Phi_1+\beta |\Phi_2|^2\Phi_1, \quad x\in \Omega, \,\,t>0,\\
-i \frac{\partial}{\partial t}\Phi_2=\Delta \Phi_2+
\mu_2|\Phi_2|^2 \Phi_2+\beta |\Phi_1|^2\Phi_2,  \quad x\in \Omega, \,\,t>0,\\
\Phi_j=\Phi_j(x,t)\in\mathbb{C},\quad j=1,2,\\
\Phi_j(x,t)=0,  \quad x\in \partial\Om, \,\,t>0, \,\,j=1,2,\end{cases}\ee
where $\Om\subset\RN (N=2, 3)$ is a smooth bounded domain, $i$ is the imaginary unit,
$\mu_1,\mu_2 >0$ and $\beta\neq 0$ is a coupling constant.
System (\ref{eq1}) arises
in mathematical models from several physical
phenomena, especially in nonlinear optics.
Physically, the solution $\Phi_j$ denotes the $j^{th}$ component of the beam
in Kerr-like photorefractive media (cf. \cite{AA}).
The positive constant $\mu_j$ is for self-focusing in the $j^{th}$ component of
the beam, and the coupling constant $\beta$ is
the interaction between the two components of the beam. Problem (\ref{eq1}) also arises
in the Hartree-Fock theory for a double
condensate, i.e., a binary mixture of Bose-Einstein
condensates in two different hyperfine states $|1\rangle$ and $|2\rangle$ (cf. \cite{EGBB}).
Physically, $\Phi_j$ are the corresponding condensate
amplitudes, $\mu_j$ and $\beta$ are the intraspecies and
interspecies scattering lengths. Precisely, the sign of $\mu_j$ represents
the self-interactions of the single state $|j\rangle$.
If $\mu_j>0$ as considered here, it is called the focusing case, in
opposition to the defocusing case where $\mu_j<0$. Besides,
the sign of $\beta$ determines whether the interactions of
states $|1\rangle$ and $|2\rangle$
are repulsive or attractive, i.e., the interaction is
attractive if $\beta>0$, and the interaction is repulsive
if $\beta<0$.

To study solitary wave solutions of (\ref{eq1}),
we set $\Phi_j(x, t)=e^{i\la_j t}u_j(x)$ for $j=1, 2$. Then system (\ref{eq1})
is reduced to the following elliptic system
\be\label{eq2}
\begin{cases}-\Delta u_1 +\la_1 u_1 =
\mu_1 u_1^3+\beta u_1 u_2^2, \quad x\in \Omega,\\
-\Delta u_2 +\la_2 u_2 =\mu_2 u_2^3+\beta u_1^2 u_2,  \quad   x\in\Om,\\
u_1=u_2=0  \,\,\,\hbox{on \,$\partial\Om$}.\end{cases}\ee

\begin{definition}\label{definition}
We call a solution $(u_1, u_2)$ {\it nontrivial}
if $u_j\not\equiv 0$ for $j=1, 2$, a solution $(u_1, u_2)$ {\it semi-trivial} if $(u_1, u_2)$ is type
of $(u_1, 0)$ or $(0, u_2)$. A solution $(u_1, u_2)$ is
called {\it positive} if $u_j > 0$ in $\Om$ for $j=1, 2$, a solution $(u_1, u_2)$
{\it sign-changing} if both $u_1$ and $u_2$ change sign,
a solution $(u_1, u_2)$ semi-nodal if one component is positive and the other one changes sign.
\end{definition}

In the last decades, system (\ref{eq2}) has received
great interest from many mathematicians.
In particular, the existence and multiplicity
of {\it positive} solutions of (\ref{eq2})
have been well studied in \cite{AC2, BDW, BW, BWW, CZ1, DWW, LW1, LW2, LWZ, MMP, NR, NTTV, S, TV, WW2} and references therein.
Note that all these papers deal with the subcritical case $N\le 3$. Recently,
Chen and Zou \cite{CZ} studied the existence and properties of
{\it positive} least energy solutions of (\ref{eq2}) in the critical case $N=4$.

On the other hand, there are few results about the existence of sign-changing or semi-nodal solutions to (\ref{eq2}) in the literature.
When $\bb>0$ is sufficiently large, multiple radially symmetric sign-changing solutions
of (\ref{eq2}) were constructed in \cite{MMP1} for the entire space case. Remark that the method in \cite{MMP1} can not be applied in the non-radial bounded domain case.
Recently, the authors \cite{CLZ} proved the existence of infinitely
many sign-changing solutions of (\ref{eq2}) for each fixed $\bb<0$. Independently, Liu, Liu and Wang \cite{LLW} obtained
infinitely many sign-changing solutions of a general $m$-coupled system ($m \ge 2$) for each fixed $\bb<0$.
The methods in \cite{CLZ, LLW} are completely different.

{\it The main goal of this paper is to study the existence of sign-changing and semi-nodal solutions when $\bb>0$ is small}.
This will complement the study made in \cite{CLZ, LLW, MMP1}. Our first result is as follows.

\bt\label{th1} Let $\Om\subset \RN (N=2, 3)$ be a smooth bounded
domain and $\la_1, \la_2$, $\mu_1, \mu_2>0$. Then for any $k\in\mathbb{N}$ there exists $\bb_k>0$ such that
system (\ref{eq2}) has at least $k$ sign-changing solutions for each fixed $\bb\in (0, \bb_k)$.\et

\begin{definition}
A nontrivial solution is called a least energy solution,
if it has the least energy among all nontrivial solutions.
A sign-changing solution is called a least energy sign-changing solution,
if it has the least energy among all sign-changing solutions.
\end{definition}

Lin and Wei \cite{LW2} proved that there exists $\bb_0>0$ small such that, for any $\bb\in (-\iy, \bb_0)$, (\ref{eq2}) has a least energy solution
which turns out to be positive. Recently, the existence of least energy sign-changing solutions for $\bb<0$ was proved in \cite{CLZ}. Here we can prove the
following result.

\bt\label{th2} Let assumptions in Theorem \ref{th1} hold. Then there exists $\bb_1'\in (0, \bb_1]$ such that system (\ref{eq2}) has a least energy sign-changing solution for each $\bb\in (0, \bb_1')$.\et

Theorems \ref{th1} and \ref{th2} are both concerned with sign-changing solutions.
The following result is about the existence of multiple semi-nodal solutions.

\bt\label{th3} Let assumptions in Theorem \ref{th1} hold. Then for any $k\in\mathbb{N}$ there exists $\bb_k>0$ such that, for each $\bb\in (0, \bb_k)$,
system (\ref{eq2}) has at least $k$ semi-nodal solutions with the first component sign-changing and the second component positive.
\et

\br\label{rmk3-1} Similarly we can prove that (\ref{eq2}) has at least $k$ semi-nodal solutions with the first component positive and the second one sign-changing for each $\bb\in (0, \bb_k)$.
Recently, \cite[Theorem 0.2]{SW} proved the existence of $\bb_k>0$ such that, for each $\bb\in (0, \bb_k)$, (\ref{eq2}) has at least $k$
nontrivial solutions $(u_{1, i}, u_{2, i})$ with $u_{1, i}>0$ in $\Om \,(i=1, \cdots, k)$.
These solutions are called semi-positive solutions in \cite{SW}. Remark that whether $u_{2, i}$ is positive or sign-changing is not known in \cite{SW}, hence
our result improves \cite[Theorem 0.2]{SW} clearly. Our proofs here
are quite different from \cite{SW}.\er

\br
Theorems \ref{th1}-\ref{th3} are all stated in the bounded domain case.
Consider
the following elliptic system in the entire space:
\be\label{eq6-1}
\begin{cases}-\Delta u_1 +\la_1 u_1 =
\mu_1 u_1^3+\beta u_1 u_2^2, \quad x\in \RN,\\
-\Delta u_2 +\la_2 u_2 =\mu_2 u_2^3+\beta u_1^2 u_2,  \quad   x\in\RN,\\
u_1(x),\,u_2(x)\to 0  \,\,\,\hbox{as}\,\,|x|\to +\iy.\end{cases}\ee
Then by working in the space $H_r^1(\RN):=\{u\in H^1(\RN) \,: \,u\,\hbox{radially symmetric}\}$ and recalling the compactness of $H_r^1(\RN)\hookrightarrow L^4(\RN)$, we can prove the following result via the same method:
For any $k\in\mathbb{N}$ there exists $\bb_k>0$
such that, for each fixed $\bb\in (0, \bb_k)$,
system (\ref{eq6-1}) has at least $k$ radially symmetric sign-changing solutions and $k$ radially symmetric semi-nodal solutions. On the other hand, in 2008 Liu and Wang \cite{LW} proved the existence of $\bb_k>0$ such that, for each $\bb\in (0, \bb_k)$, (\ref{eq6-1}) has at least $k$ nontrivial radially symmetric solutions. In fact, they studied a general $m$-coupled system ($m\ge 2$). Remark that whether solutions obtained in \cite{LW} are positive or sign-changing or semi-nodal is not known. Moreover, Liu and Wang \cite[Remark 3.6]{LW} suspected that solutions obtained in \cite{LW} are not positive solutions, but no proof has yet been given. Our results improve the result of \cite{LW} in the two coupled case (m=2).
\er

\br After the completion of this paper, we learned that (\ref{eq6-1}) has also been
studied in a recent manuscript \cite{KKL}, where the authors obtained multiple radially symmetric sign-changing solutions with a prescribed number of zeros for $\bb>0$ small. Remark that their method can not be applied in the non-radial bounded domain case.\er

The rest of this paper proves these theorems. We give some notations here. Throughout this paper,
we denote the norm of $L^p(\Om)$ by $|u|_p =
(\int_{\Om}|u|^p\,dx)^{\frac{1}{p}}$, the norm of $H^1_0(\Om)$ by $\|u\|^2=\int_{\Om}(|\nabla u|^2+u^2)\,dx$ and positive constants
(possibly different in different places) by $C, C_0, C_1, \cdots$. Denote
$$\|u\|_{\la_i}^2:=\int_{\Om} (|\nabla u|^2 + \la_i u^2)\,dx$$
for convenience. Then $\|\cdot\|_{\la_i}$ are
equivalent norms to $\|\cdot\|$. Define $H:=H_0^1(\Om)\times H_0^1(\Om)$ with norm $\|(u_1, u_2)\|_H^2:=\|u_1\|_{\la_1}^2+\|u_2\|_{\la_2}^2$.

The rest of this paper is organized as follows. In Section 2 we give the proof of Theorem \ref{th1}.
The main ideas of our proof are inspired by \cite{CLZ, TT},
where a new {\it constrained problem} introduced by \cite{CLZ} and a new
notion of {\it vector genus} introduced by \cite{TT} will be used to define appropriate minimax values.
In \cite{TT}, Tavares and Terracini
studied the following general $m$-coupled system
\be\label{eq3}
\begin{cases}-\Delta u_j -\mu_j u_j^3-\bb u_j\sum_{i\neq j} u_i^2 =\la_{j, \bb} u_j,\\
u_j\in H_0^1(\Om),\quad j=1,\cdots, m,\end{cases}\ee
where $\bb<0$, $\mu_j\le 0$ are all fixed constants. Then \cite[Theorem 1.1]{TT} says that there exist
infinitely many $\la=(\la_{1, \bb},\cdots,\la_{m, \bb})\in \R^m$ and $u=(u_1, \cdots, u_m)\in H_0^1(\Om, \R^m)$
such that $(u, \la)$ are sign-changing solutions of (\ref{eq3}).
That is, $\la_{j, \bb}$ is not
fixed a {\it priori} and appears as a Lagrange multiplier in \cite{TT}.
Here we deal with the focusing case $\mu_j>0$,
and $\la_j, \mu_j$, $\bb>0$ are all fixed constants.
Some arguments in our proof are borrowed from \cite{CLZ, TT} with modifications.
Although some procedures are close to
those in \cite{CLZ, TT}, we prefer to provide all the necessary details to make the paper
self-contained. In Section 3 we will use a minimizing argument to prove Theorem \ref{th2}.
By giving some modifications to arguments in Sections 2 and 3, we will prove
Theorems \ref{th3} in Section 4.

\vskip0.1in

\s{Proof of Theorem \ref{th1}}
\renewcommand{\theequation}{2.\arabic{equation}}

In the sequel we let assumptions in Theorem \ref{th1} hold.
Without loss of generality we assume that $\mu_1\ge \mu_2$.
Let $\bb\in (0, \mu_2)$. Note that solutions
of (\ref{eq2}) correspond to the critical points of $C^2$ functional $E_\bb: H\to \R$
given by
{\allowdisplaybreaks
\begin{align}\label{eq2-1}
E_\bb(u_1, u_2):=\frac{1}{2} \left(\|u_1\|_{\la_1}^2+\|u_2\|_{\la_2}^2\right)-\frac{1}{4}\left(\mu_1 |u_1|_4^4+\mu_2|u_2|_4^4\right)-\frac{\beta}{2} \int_{\Om}u_1^2u_2^2\,dx.
\end{align}
}%
Since we are only concerned with nontrivial solutions,
we denote $\wH:=\{(u_1, u_2)\in H : u_i\neq 0\,\,\hbox{for}\,\, i=1, 2\}$, which is open in $H$. Write $\vec{u}=(u_1, u_2)$ for convenience.

\bl\label{lemma1}For any $\vec{u}=(u_1, u_2)\in\wH$, if
\be\label{eq2-2}
\begin{cases}
\mu_2|u_2|_4^4\|u_1\|_{\la_1}^2-\bb\|u_2\|_{\la_2}^2\int_{\Om}u_1^2u_2^2\,dx>0,\\
\mu_1|u_1|_4^4\|u_2\|_{\la_2}^2-\bb\|u_1\|_{\la_1}^2\int_{\Om}u_1^2u_2^2\,dx>0,
\end{cases}\ee
then system
\be\label{eq2-3}
\begin{cases}\|u_1\|_{\la_1}^2 =
t_1\mu_1 |u_1|_4^4+t_2\beta\int_{\Om} u_1^2 u_2^2\,dx\\
\|u_2\|_{\la_2}^2 =
t_2\mu_2 |u_2|_4^4+t_1\beta\int_{\Om} u_1^2 u_2^2\,dx\end{cases}\ee
has a unique solution
\be\label{eq2-4}
\begin{cases}t_1(\vec{u})=\frac{\mu_2|u_2|_4^4\|u_1\|_{\la_1}^2-\bb\|u_2\|_{\la_2}^2\int_{\Om}u_1^2u_2^2\,dx}{\mu_1\mu_2|u_1|_4^4 |u_2|_4^4-\bb^2(\int_{\Om}u_1^2 u_2^2\,dx)^2}>0\\
t_2(\vec{u})=\frac{\mu_1|u_1|_4^4\|u_2\|_{\la_2}^2-\bb\|u_1\|_{\la_1}^2\int_{\Om}u_1^2u_2^2\,dx}{\mu_1\mu_2|u_1|_4^4 |u_2|_4^4-\bb^2(\int_{\Om}u_1^2 u_2^2\,dx)^2}>0.\end{cases}\ee
Moreover,
{\allowdisplaybreaks
\begin{align}\label{eq2-5}\sup_{t_1, t_2\ge 0}&E_\bb \left(\sqrt{t_1}u_1, \,\sqrt{t_2} u_2\right)=E_\bb\left(\sqrt{t_1(\vec{u})}u_1, \sqrt{t_2(\vec{u})}u_2\right)\nonumber\\
&=\frac{1}{4}\left(t_1(\vec{u})\|u_1\|_{\la_1}^2+t_2(\vec{u})\|u_2\|_{\la_2}^2\right)\nonumber\\
&=\frac{1}{4}\frac{\mu_2|u_2|_4^4\|u_1\|_{\la_1}^4
-2\bb\|u_1\|_{\la_1}^2\|u_2\|_{\la_2}^2\int_{\Om}u_1^2u_2^2\,dx+\mu_1|u_1|_4^4\|u_2\|_{\la_2}^4}
{\mu_1\mu_2|u_1|_4^4 |u_2|_4^4-\bb^2(\int_{\Om}u_1^2 u_2^2\,dx)^2}
\end{align}
}%
and $(t_1(\vec{u}), t_2(\vec{u}))$ is the unique maximum point of $E_\bb (\sqrt{t_1}u_1, \sqrt{t_2}u_2)$.
\el

\noindent {\bf Proof. } By (\ref{eq2-2}) we see that
$\mu_1\mu_2|u_1|_4^4 |u_2|_4^4-\bb^2(\int_{\Om}u_1^2 u_2^2dx)^2>0$,
so $(t_1(\vec{u}), t_2(\vec{u}))$ defined in (\ref{eq2-4}) is the unique solution of (\ref{eq2-3}).
Note that for $t_1, t_2\ge 0$,
{\allowdisplaybreaks
\begin{align*}
f(t_1, t_2):=& E_\bb\left(\sqrt{t_1}u_1, \sqrt{t_2}u_2\right)=\frac{1}{2} t_1\|u_1\|_{\la_1}^2+
\frac{1}{2}t_2\|u_2\|_{\la_2}^2\\
&-\frac{1}{4}\left(t_1^2\mu_1 |u_1|_4^4+ t_2^2\mu_2|u_2|_4^4\right)-\frac{1}{2}t_1 t_2\beta \int_{\Om}u_1^2u_2^2\,dx\\
\le & \left(\frac{t_1}{2}\|u_1\|_{\la_1}^2-\frac{t_1^2}{4}\mu_1 |u_1|_4^4\right)+
\left(\frac{t_2}{2}\|u_2\|_{\la_2}^2-\frac{t_2^2}{4}\mu_2 |u_2|_4^4\right).
\end{align*}
}%
This implies that $f (t_1, t_2)< 0$ for $\max\{t_1, t_2\}> T$, where $T$ is some positive constant.
So there exists $(\tilde{t}_1, \tilde{t}_2)\in [0, T]^2\setminus\{(0, 0)\}$ such that
$$f(\tilde{t}_1, \tilde{t}_2)=\sup_{t_1, t_2\ge 0}f (t_1, t_2).$$
It suffices to show that $(\tilde{t}_1, \tilde{t}_2)=(t_1(\vec{u}), t_2(\vec{u}))$. Note that
$$\sup_{t_1\ge 0}f(t_1, 0)=\frac{1}{4}\frac{\|u_1\|_{\la_1}^4}{\mu_1 |u_1|_4^4}.$$
Recalling the expression of $f(t_1(\vec{u}), t_2(\vec{u}))$ in (\ref{eq2-5}), by a direct computation we deduce from (\ref{eq2-2}) that
$$f(t_1(\vec{u}), t_2(\vec{u}))-\sup_{t_1\ge 0}f(t_1, 0)=\frac{(\mu_1|u_1|_4^4\|u_2\|_{\la_2}^2-\bb\|u_1\|_{\la_1}^2\int_{\Om}u_1^2u_2^2\,dx)^2}{4\mu_1|u_1|_4^4[\mu_1\mu_2|u_1|_4^4 |u_2|_4^4-\bb^2(\int_{\Om}u_1^2 u_2^2\,dx)^2]}>0.$$
Similarly we have $f(t_1(\vec{u}), t_2(\vec{u}))-\sup_{t_2\ge 0}f(0, t_2)>0$, so $\tilde{t}_1>0$ and $\tilde{t}_2>0$. Then by
$\frac{\partial}{\partial t_1}f(t_1, t_2)|_{(\tilde{t}_1, \tilde{t}_2)}=\frac{\partial}{\partial t_2}f(t_1, t_2)|_{(\tilde{t}_1, \tilde{t}_2)}=0$
we see that $(\tilde{t}_1, \tilde{t}_2)$ satisfies (\ref{eq2-3}), so $(\tilde{t}_1, \tilde{t}_2)=(t_1(\vec{u}), t_2(\vec{u}))$.\hfill$\square$\\

Define
{\allowdisplaybreaks
\begin{align}\label{eq2-8-1}
&\mathcal{M}^\ast:=\left\{\vec{u}\in H\,\,:\,\, 1/2<|u_1|^4_4<2,\,\, 1/2<|u_2|^4_4<2\right\};\\
&\mathcal{M}_\bb^\ast:=\left\{\vec{u}\in \mathcal{M}^\ast\,\,:\,\, \hbox{$\vec{u}$ satisfies (\ref{eq2-2})}\right\};\nonumber\\
&\mathcal{M}_\bb^{\ast\ast}:=\left\{\vec{u}\in \mathcal{M}^\ast\,\,:\,\, \begin{array}{ll}
\mu_2\|u_1\|_{\la_1}^2-\bb\|u_2\|_{\la_2}^2\int_{\Om}u_1^2u_2^2\,dx>0\\
\mu_1\|u_2\|_{\la_2}^2-\bb\|u_1\|_{\la_1}^2\int_{\Om}u_1^2u_2^2\,dx>0
\end{array}\right\};\nonumber\\
\label{eq2-8}&\mathcal{M}:=\left\{\vec{u}\in H\,\,:\,\, |u_1|_4=1,\,\, |u_2|_4=1\right\},\quad \mathcal{M}_\bb:=\mathcal{M}\cap \mathcal{M}_\bb^\ast.
\end{align}
}%
Then $\mathcal{M}_\bb=\mathcal{M}\cap \mathcal{M}_\bb^{\ast\ast}$.
Evidently $\mathcal{M}^\ast$, $\mathcal{M}_\bb^*$, $\mathcal{M}_\bb^{\ast\ast}$ are all open subsets of $H$ and $\mathcal{M}$ is closed.
Note that $\mu_1\mu_2-\bb^2(\int_{\Om}u_1^2 u_2^2\,dx)^2>0$ for any $\vec{u}\in\mathcal{M}_\bb^{\ast\ast}$, as in \cite{CLZ} we
define a new functional $J_\bb : \mathcal{M}_\bb^{\ast\ast}\to (0, +\iy)$ by
$$J_\bb (\vec{u}):=\frac{1}{4}\frac{\mu_2\|u_1\|_{\la_1}^4
-2\bb\|u_1\|_{\la_1}^2\|u_2\|_{\la_2}^2\int_{\Om}u_1^2u_2^2\,dx+\mu_1\|u_2\|_{\la_2}^4}{\mu_1\mu_2-\bb^2(\int_{\Om}u_1^2 u_2^2\,dx)^2}.
$$
A direct computation yields
$J_\bb\in C^1(\mathcal{M}^{\ast\ast}_\bb, \,(0, +\iy))$. Moreover, since any $\vec{u}\in\mathcal{M}_\bb$ is an interior point of $\mathcal{M}^{\ast\ast}_\bb$, by (\ref{eq2-4}) we can prove that
{\allowdisplaybreaks
\begin{align}
\label{eq2-11}&J_\bb'(\vec{u})(\vp, 0)=t_1(\vec{u})\int_{\Om}(\nabla u_1\nabla\vp+\la_1 u_1\vp)\,dx-t_1(\vec{u})t_2(\vec{u})\bb\int_{\Om}u_1 u_2^2\vp\,dx,\\
\label{eq2-11-1}&J_\bb'(\vec{u})(0,\psi)=t_2(\vec{u})\int_{\Om}(\nabla u_2\nabla\psi+\la_2 u_2\psi)\,dx-t_1(\vec{u})t_2(\vec{u})\bb\int_{\Om}u_1^2 u_2\psi\,dx
\end{align}
}%
hold for any $\vec{u}\in \mathcal{M}_\bb$ and $\vp,\, \psi\in H_0^1(\Om)$
(Remark that (\ref{eq2-11})-(\ref{eq2-11-1}) do not hold for $\vec{u}\in\mathcal{M}_\bb^{\ast\ast}\setminus\mathcal{M}_\bb$). Note that
Lemma \ref{lemma1} yields
\be\label{eq2-9}J_\bb (u_1, u_2)=\sup_{t_1, t_2\ge 0}E_\bb \left(\sqrt{t_1}u_1, \,\sqrt{t_2} u_2\right),\quad \forall\,(u_1, u_2)\in \mathcal{M}_\bb.\ee
To obtain nontrivial solutions of (\ref{eq2}),
we turn to study the functional $J_\bb$ restricted to $\mathcal{M}_\bb$, which is a problem with two constraints. Define
\begin{align}\label{eq22-1}
\mathcal{N}_b^{\ast}:=\left\{\vec{u}\in\mathcal{M}^\ast \,\,:\,\, \|u_1\|_{\la_1}^2,\,\|u_2\|_{\la_2}^2<b\right\},
\quad\mathcal{N}_b:=\mathcal{N}_b^\ast\cap \mathcal{M}.
\end{align}

Fix any $k\in \mathbb{N}$. Our goal is to prove the existence of $\bb_k>0$ such
that (\ref{eq2}) has at least $k$ sign-changing solutions for any $\bb\in (0, \bb_k)$. To do this, we
let $W_{k+1}$ be a $k+1$ dimensional subspace of $H_0^1(\Om)$ which contains an element $\vp_0$ satisfying $\vp_0>0$ in $\Om$.
Then we can find $\bar{b}>0$ such that
\be\label{eq22-2}\|u\|_{\la_1}^2,\,\|u\|_{\la_2}^2<\bar{b},\quad \forall\,u\in W_{k+1}\,\,\,\hbox{satisfying}\,\,|u|^4_4<2.\ee
Fix a $b>0$ such that
\be\label{eq22-3}b^2>(2+\mu_1/\mu_2)\bar{b}^2.\ee
Then $\mathcal{N}_{\bar{b}}^{\ast}\subset\mathcal{N}_b^{\ast}$ and $\mathcal{N}_{\bar{b}}\subset\mathcal{N}_b$. Recalling the Sobolev inequality
\be\label{eq2-10}\|u\|_{\la_i}^2\ge S |u|_4^2,\quad\forall\, u\in H_0^1(\Om),\,\,\,i=1, 2,\ee
where $S$ is a positive constant, we have the following lemma.

\bl\label{lemma2}There exist $\bb_0\in (0, \mu_2)$ and $C_1>C_0>0$ such that for any $\bb\in (0, \bb_0)$ there hold
$\mathcal{N}_b^\ast\subset\mathcal{M}_\bb^{\ast}\cap\mathcal{M}_\bb^{\ast\ast}$ and
$$C_0\le t_1(\vec{u}),\, t_2(\vec{u})\le C_1,\quad\forall\, \vec{u}\in \mathcal{N}_b^\ast.$$\el

\noindent {\bf Proof. }Define
$\bb_0:=\frac{\mu_2 S}{8b}$ and let $\bb\in (0, \bb_0)$. For any $\vec{u}=(u_1, u_2)\in\mathcal{N}_b^\ast$, we see from
(\ref{eq2-8-1}) and (\ref{eq2-10}) that $\int_{\Om}u_1^2u_2^2\,dx\le |u_1|_4^2|u_2|_4^2<2$ and $\|u_i\|_{\la_i}^2\ge S/\sqrt{2}$. Hence
{\allowdisplaybreaks
\begin{align*}
&\mu_2|u_2|_4^4\|u_1\|_{\la_1}^2-\bb\|u_2\|_{\la_2}^2\int_{\Om}u_1^2u_2^2\,dx\ge\frac{\mu_2 S}{2\sqrt{2}}-2b\bb_0\ge\frac{\mu_2 S}{16};\\
&\mu_1|u_1|_4^4\|u_2\|_{\la_2}^2-\bb\|u_1\|_{\la_1}^2\int_{\Om}u_1^2u_2^2\,dx\ge\frac{\mu_2 S}{16};\\
&\mu_2\|u_1\|_{\la_1}^2-\bb\|u_2\|_{\la_2}^2\int_{\Om}u_1^2u_2^2\,dx\ge\frac{\mu_2 S}{16};\\
&\mu_1\|u_2\|_{\la_2}^2-\bb\|u_1\|_{\la_1}^2\int_{\Om}u_1^2u_2^2\,dx\ge\frac{\mu_2 S}{16};\\
&\mu_1\mu_2-\bb^2\left(\int_{\Om}u_1^2 u_2^2\,dx\right)^2\ge\frac{\mu_2^2 S^2}{2^{8}}\cdot\frac{1}{\|u_1\|_{\la_1}^2\|u_2\|_{\la_2}^2}\ge\frac{\mu_2^2 S^2}{2^{8}b^2};\\
&\mu_1\mu_2|u_1|_4^4 |u_2|_4^4-\bb^2\left(\int_{\Om}u_1^2 u_2^2\,dx\right)^2\ge\frac{\mu_2^2 S^2}{2^{8}b^2}.
\end{align*}
}%
Then $\vec{u}\in\mathcal{M}_\bb^{\ast}\cap\mathcal{M}_\bb^{\ast\ast}$. Moreover, combining these with (\ref{eq2-4}) we have
$$t_i(\vec{u})\ge \frac{\mu_2 S}{2^4}\cdot\frac{1}{\mu_1\mu_2|u_1|_4^4|u_2|_4^4}\ge\frac{S}{2^{6}\mu_1},
\,\,\, t_i(\vec{u})\le \frac{2^{9}b^3}{\mu_2^2S^2}\mu_1,\,\,\, i=1, 2.$$
This completes the proof.\hfill$\square$

\bl\label{lemma22}There exist $\bb_k\in (0, \bb_0]$ and $d_k>0$ such that
\be\label{eq22-5}\inf_{\vec{u}\in \partial\mathcal{N}_b}J_\bb(\vec{u})\ge d_k>\sup_{\vec{u}\in\mathcal{N}_{\bar{b}}} J_\bb(\vec{u}),\quad\forall\,\bb\in(0, \bb_k).
\ee\el

\noindent {\bf Proof. }This proof is inspired by \cite{SW}. Define
$$I_i(u_i):=\frac{1}{4\mu_i}\|u_i\|_{\la_i}^4,\quad i=1, 2.$$
Then for any $\vec{u}\in\overline{\mathcal{N}_b}$ and $\bb\in (0, \bb_0)$ we have
{\allowdisplaybreaks
\begin{align*}
&|J_\bb(\vec{u})-I_1(u_1)-I_2(u_2)|\\
=&\frac{\bb\left|\bb(\int_{\Om}u_1^2 u_2^2\,dx)^2\sum_{i=1}^2\|u_i\|_{\la_i}^4/\mu_i-2\|u_1\|_{\la_1}^2\|u_2\|_{\la_2}^2\int_{\Om}u_1^2 u_2^2\,dx\right|}{4[\mu_1\mu_2-\bb^2(\int_{\Om}u_1^2 u_2^2\,dx)^2]}
\le  C\bb,
\end{align*}
}%
where $C>0$ is independent of $\vec{u}\in\overline{\mathcal{N}_b}$ and $\bb\in (0, \bb_0)$. Therefore,
\begin{align*}&\sup_{\vec{u}\in\mathcal{N}_{\bar{b}}} J_\bb(\vec{u})\le\sup_{\vec{u}\in\mathcal{N}_{\bar{b}}}(I_1(u_1)+I_2(u_2))+C\bb
\le\frac{\bar{b}^2}{4\mu_1}+\frac{\bar{b}^2}{4\mu_2}+C\bb;\\
&\inf_{\vec{u}\in \partial\mathcal{N}_b}J_\bb(\vec{u})\ge\inf_{\vec{u}\in \partial\mathcal{N}_b}(I_1(u_1)+I_2(u_2))-C\bb\ge \frac{b^2}{4\mu_1}-C\bb.\end{align*}
Recalling (\ref{eq22-3}), we let $\bb_k=\min\{\frac{\bar{b}^2}{8\mu_1 C}, \bb_0\}$ and $d_k=\frac{b^2}{4\mu_1}-C\bb_k$, then (\ref{eq22-5}) holds. This completes the proof.
\hfill$\square$\\

In the following we always let $(i, j)=(1, 2)$ or $(i, j)=(2, 1)$. Recalling (\ref{eq2-10}) and Lemma \ref{lemma2},
we can take $\bb_k$ smaller if necessary such that, for any $\bb\in (0, \bb_k)$ and $\vec{u}\in\mathcal{N}_b^\ast$, there holds
\be\label{eq22-4}\|v\|_{\la_i}^2-\bb t_j(\vec{u})\int_{\Om}u_j^2v^2\,dx\ge \frac{1}{2}\|v\|_{\la_i}^2,\quad\forall\,v\in H_0^1(\Om),\,\,i=1,2.\ee
Clearly (\ref{eq22-4}) implies that the operators $-\Delta+\la_i-\bb t_j(\vec{u})u_j^2$ are positive definite in $H_0^1(\Om)$.
In the rest of this section we fix any $\bb\in (0, \bb_k)$. We will show that (\ref{eq2}) has at least $k$ sign-changing solutions.
For any $\vec{u}=(u_1, u_2)\in \mathcal{N}_b^\ast$, let $\tilde{w}_i\in H_0^1(\Om)$ be the unique
solution of the following linear problem
\be\label{eq2-012}-\Delta \tilde{w}_i+\la_i \tilde{w}_i-\bb t_j(\vec{u})u_j^2 \tilde{w}_i=\mu_i t_i(\vec{u})u_i^3,\quad \tilde{w}_i\in H_0^1(\Om).\ee
Since $|u_i|^4_4 > 1/2$, so $\tilde{w}_i\neq 0$ and we see from (\ref{eq22-4}) that
$$\int_{\Om}u_i^3 \tilde{w}_i\,dx=\frac{1}{\mu_i t_i(\vec{u})}\left(\|\tilde{w}_i\|_{\la_i}^2-\bb t_j(\vec{u})\int_{\Om}u_j^2 \tilde{w}_i^2\,dx\right)
\ge \frac{1}{2\mu_i t_i(\vec{u})}\|\tilde{w}_i\|_{\la_i}^2>0.$$
Define
\be\label{eq2-12}w_i=\al_i \tilde{w}_i,\quad\hbox{where}\,\,\al_i=\frac{1}{\int_{\Om}u_i^3 \tilde{w}_i\,dx}>0.\ee
Then $w_i$ is the unique solution of the following problem
\be\label{eq2-13}\begin{cases}
-\Delta w_i+\la_i w_i-\bb t_j(\vec{u})u_j^2 w_i=\al_i\mu_i t_i(\vec{u})u_i^3,\quad w_i\in H_0^1(\Om),\\
\int_{\Om}u_i^3 w_i\,dx=1.
\end{cases}\ee
Now we define an operator $K=(K_1, K_2) : \mathcal{N}_b^\ast\to H$ by
\be\label{eq2-14}K(\vec{u})=(K_1(\vec{u}), K_2(\vec{u})):=\vec{w}=(w_1, w_2).\ee
Define the transformations
\be\label{involution}\sg_i: H\to H\quad\hbox{by}\quad \sg_1(u_1, u_2):=(-u_1, u_2),\,\,\,\sg_2(u_1, u_2):=(u_1, -u_2).\ee
Then it is easy to check that
\be\label{eq2-15}K(\sg_i(\vec{u}))=\sg_i(K(\vec{u})),\quad i=1, 2.\ee

\bl\label{lemma3} $K\in C^1(\mathcal{N}_b^\ast, H)$.\el

\noindent {\bf Proof. } It suffices to apply the Implicit Theorem to the $C^1$ map
{\allowdisplaybreaks
\begin{align*}
&\Psi : \mathcal{N}_b^\ast\times H_0^1(\Om)\times \R\to H_0^1(\Om)\times\R, \quad\hbox{where}\\
&\Psi(\vec{u}, v, \al)=\left(v-(-\Delta+\la_i)^{-1}\left(\bb t_j(\vec{u})u_j^2 v+\al \mu_i t_i(\vec{u})u_i^3\right),\,\,\int_{\Om}u_i^3v\,dx -1\right).
\end{align*}
}%
Note that (\ref{eq2-13}) holds if and only if $\Psi(\vec{u}, w_i, \al_i)=(0, 0)$. By computing the derivative of $\Psi$ with respect to $(v, \al)$ at the
point $(\vec{u}, w_i, \al_i)$ in the direction $(\bar{w}, \bar{\al})$, we obtain a map $\Phi: H_0^1(\Om)\times \R\to H_0^1(\Om)\times\R$ given by
{\allowdisplaybreaks
\begin{align*}
\Phi(\bar{w}, \bar{\al}):=& D_{v, \al}\Psi(\vec{u}, w_i, \al_i)(\bar{w}, \bar{\al})\\
=& \left(\bar{w}-(-\Delta+\la_i)^{-1}
\left(\bb t_j(\vec{u})u_j^2 \bar{w}+\bar{\al} \mu_i t_i(\vec{u})u_i^3\right),\,\,\int_{\Om}u_i^3\bar{w}\,dx\right).
\end{align*}
}%
Recalling (\ref{eq22-4}), similarly as \cite[Lemma 2.3]{CLZ} it is easy to prove that $\Phi$ is a bijective map. We omit the details.\hfill$\square$

\bl\label{lemma4}Assume that $\{\vec{u}_n=(u_{n, 1}, u_{n, 2}) : n\ge 1\}\subset \mathcal{N}_b$.
Then there exists $\vec{w}\in H$ such that, up to a subsequence,
$\vec{w}_n:=K(\vec{u}_n)\to \vec{w}$ strongly in $H$.\el

\noindent {\bf Proof. }
Up to a subsequence, we may assume that $\vec{u}_n\rightharpoonup \vec{u}=(u_1, u_2)$ weakly in $H$ and so
$u_{n, i}\to u_i$ strongly in $L^{4}(\Om)$, which implies $|u_i|_4=1$.  Moreover, by Lemma \ref{lemma2} we may assume $t_i(\vec{u}_n)\to t_i>0$.
Recall that $w_{n, i}=\al_{n, i}\tilde{w}_{n, i}$, where $\al_{n, i}$ and $\tilde{w}_{n, i}$ are seen in (\ref{eq2-012})-(\ref{eq2-12}).
By (\ref{eq22-4})-(\ref{eq2-012}) we have
$$\frac{1}{2}\|\tilde{w}_{n, i}\|_{\la_i}^2\le \mu_i t_i(\vec{u}_n)\int_{\Om}u_{n, i}^3 \tilde{w}_{n, i}\,dx\le C|\tilde{w}_{n, i}|_4\le C\|\tilde{w}_{n, i}\|_{\la_i},$$
which implies that $\{\tilde{w}_{n, i} : n\ge 1\}$ are bounded in $H_0^1(\Om)$. Up to a subsequence, we may assume that
$\tilde{w}_{n, i}\to \tilde{w}_i$ weakly in $H_0^1(\Om)$ and strongly in $L^4(\Om)$.
Then by (\ref{eq2-012}) it is standard to prove that $\tilde{w}_{n, i}\to \tilde{w}_i$ strongly in $H_0^1(\Om)$. Moreover, $\tilde{w}_i$ satisfies
$-\Delta \tilde{w}_i+\la_i \tilde{w}_i-\bb t_ju_j^2 \tilde{w}_i=\mu_i t_i u_i^3$.
Since $|u_i|_4=1$, so $\tilde{w}_i\neq 0$ and then $\int_{\Om}u_i^3 \tilde{w}_i\,dx>0$, which implies that
$$\lim_{n\to\iy}\al_{n, i}=\lim_{n\to\iy}\frac{1}{\int_{\Om}u_{n,i}^3 \tilde{w}_{n,i}\,dx}=\frac{1}{\int_{\Om}u_i^3 \tilde{w}_i\,dx}=:\al_i.$$
Therefore, $w_{n, i}=\al_{n, i}\tilde{w}_{n, i}\to \al_i\tilde{w_i}=: w_i$ strongly in $H_0^1(\Om)$.\hfill$\square$\\

To continue our proof, we need to use {\it vector genus} introduced by \cite{TT} to define proper minimax energy levels.
Recall (\ref{eq2-8}) and  (\ref{involution}), as in \cite{TT} we consider the class of sets
$$\mathcal{F}=\{A\subset \mathcal{M} : A\,\,\hbox{ is closed and}\,\,\sg_i(\vec{u})\in A\,\,\forall\,\vec{u}\in A,\,\,i=1, 2\},$$
and, for each $A\in \mathcal{F}$ and $k_1, k_2\in\mathbb{N}$, the class of functions
$$F_{(k_1, k_2)}(A)=\left\{f=(f_1, f_2): A\to \prod_{i=1}^2\R^{k_i-1} : \begin{array}{lll}  f_i : A\to \R^{k_i-1}\,\,\hbox{continuous,}\\
 f_i(\sg_i(\vec{u}))=-f_i(\vec{u})\,\,\hbox{for each}\,\,i,\\
 f_i(\sg_j(\vec{u}))=f_i(\vec{u})\,\,\hbox{for}\,\,j\neq i
\end{array} \right\}.$$
Here we denote $\R^0:=\{0\}$. Let us recall vector genus from \cite{TT}.
\begin{definition}\label{definition1} (Vector genus, see \cite{TT})
Let $A\in \mathcal{F}$ and take any $k_1, k_2\in\mathbb{N}$. We say that $\vec{\ga}(A)\ge (k_1, k_2)$ if for every $f\in F_{(k_1, k_2)}(A)$ there exists $\vec{u}\in A$ such that $f(\vec{u})=(f_1(\vec{u}), f_2(\vec{u}))=(0, 0)$. We denote
$$\Gamma^{(k_1, k_2)}:=\{A\in\mathcal{F} : \vec{\ga}(A)\ge (k_1, k_2)\}.$$
\end{definition}

\bl\label{lemma5}(see \cite{TT}) With the previous notations, the following properties hold.
\begin{itemize}
\item[$(i)$] Take $A_1\times A_2\subset \mathcal{M}$ and let $\eta_i: S^{k_i-1}:=\{x\in \R^{k_i} : |x|=1\}\to A_i$ be a homeomorphism such that
$\eta_i(-x)=-\eta_i(x)$ for every $x\in S^{k_i-1}$, $i=1, 2$. Then $A_1\times A_2\in \Gamma^{(k_1, k_2)}$.

\item[$(ii)$] We have $\overline{\eta(A)}\in \Gamma^{(k_1, k_2)}$ whenever $A\in \Gamma^{(k_1, k_2)}$ and a continuous map $\eta : A\to \mathcal{M}$ is such
that $\eta\circ \sg_i=\sg_i\circ\eta,\,\,\forall\,i=1, 2$.
\end{itemize}\el

To obtain sign-changing solutions, as in many references such as \cite{CMT, BLT, Zou}, we should use cones of positive functions.
Precisely, we define
\be\label{cone}\mathcal{P}_i:=\{\vec{u}=(u_1, u_2)\in H : u_i\ge 0\},\quad \mathcal{P}:=\bigcup_{i=1}^2 (\mathcal{P}_i\cup -\mathcal{P}_i).\ee
Moreover, for $\dd>0$ we define
$\mathcal{P}_\dd:=\{\vec{u}\in H \,:\,\hbox{dist}_4(\vec{u}, \mathcal{P})<\dd \}$,
where
{\allowdisplaybreaks
\begin{align}\label{cone1}&\hbox{dist}_4(\vec{u}, \mathcal{P}):=\min\big\{\hbox{dist}_4(u_i,\,\mathcal{P}_i),\,\,\hbox{dist}_4(u_i,\,-\mathcal{P}_i),\quad i=1, 2\big\},\\
&\hbox{dist}_4(u_i,\,\pm\mathcal{P}_i ):=\inf\{|u_i-v|_4 \,\,:\,\, v\in \pm\mathcal{P}_i\}.\nonumber
\end{align}
}%
Denote $u^{\pm}:=\max\{0, \pm u\}$, then it is easy to check that $\hbox{dist}_4(u_i, \pm\mathcal{P}_i)=|u_i^{\mp}|_4$.
The following lemma was proved in \cite{CLZ}.

\bl\label{lemma6} (see \cite[Lemma 2.6]{CLZ}) Let $k_1, k_2\ge 2$. Then for any $\dd<2^{-1/4}$ and any $A\in\Gamma^{(k_1, k_2)}$ there holds $A\setminus \mathcal{P}_\dd\neq\emptyset$.\el

\bl\label{lemma7}  There exists $A\in \Gamma^{(k+1, k+1)}$ such that $A\subset \mathcal{N}_b$ and
$\sup_A J_\bb< d_k$. \el

\noindent {\bf Proof. } Recalling $W_{k+1}$ in (\ref{eq22-2}), we define
$$A_1=A_2:=\big\{u\in W_{k+1} \,:\, |u|_4=1\big\}.$$
Note that there exists an obvious odd homeomorphism from $S^{k}$ to $A_i$. By Lemma \ref{lemma5}-$(i)$ one has $A:=A_1\times A_2\in \Gamma^{(k+1, k+1)}$. We see from (\ref{eq22-2}) that $A\subset \mathcal{N}_{\bar{b}}$, and so Lemma \ref{lemma22} yields $\sup_{A}J_\bb < d_k$.
\hfill$\square$\\

For every $k_1, k_2\in [2, k+1]$ and $0<\dd < 2^{-1/4}$, we define
\be\label{eq2-17}c_{\bb,\dd}^{k_1, k_2}:=\inf_{A\in \Gamma_\bb^{(k_1, k_2)}}\sup_{\vec{u}\in A\setminus \mathcal{P}_\dd}J_\bb(\vec{u}),\ee
where
\be\label{eq2-17-1}\Gamma_\bb^{(k_1, k_2)}:=\left\{A\in \Gamma^{(k_1, k_2)} \,:\, A\subset \mathcal{N}_b,\,\,\,\sup_{A}J_\bb< d_k\right\}.\ee
Noting that $\Gamma_\bb^{(\tilde{k}_1, \tilde{k}_2)}\subset \Gamma_\bb^{(k_1, k_2)}$ for any $\tilde{k}_1\ge k_1$ and $\tilde{k}_2\ge k_2$, we see that Lemma \ref{lemma7} yields $\Gamma_\bb^{(k_1, k_2)}\neq \emptyset$ and so $c_{\bb, \dd}^{k_1, k_2}$ is well defined for any $k_1, k_2\in [2, k+1]$. Moreover,
$$c_{\bb, \dd}^{k_1, k_2}< d_k\quad\hbox{for every}\,\,\dd\in (0, 2^{-1/4})\,\,\hbox{and}\,\,k_1, k_2\in [2, k+1]. $$
We will prove that $c_{\bb, \dd}^{k_1, k_2}$ is a critical value of $E_\bb$ for $\dd>0$ sufficiently small. Define $\mathcal{N}_{b,\bb}:=\{\vec{u}\in  \mathcal{N}_b: J_\bb(\vec{u})<d_k\}$,
then Lemma \ref{lemma22} yields $\mathcal{N}_{\bar{b}}\subset\mathcal{N}_{b, \bb}$.

\bl\label{lemma8}For any sufficiently small $\dd\in (0, 2^{-1/4})$, there holds
$$\hbox{dist}_4(K(\vec{u}), \mathcal{P})<\dd/2,\quad\forall\,\vec{u}\in\mathcal{N}_{b, \bb},\,\,
\hbox{dist}_4(\vec{u}, \mathcal{P})<\dd.$$\el

\noindent {\bf Proof. } Assume by contradiction that there exist $\dd_n\to 0$ and $\vec{u}_n=(u_{n, 1}, u_{n, 2})\in\mathcal{N}_{b,\bb}$
such that
$\hbox{dist}_4(\vec{u}_n, \mathcal{P})<\dd_n$ and $\hbox{dist}_4(K(\vec{u}_n), \mathcal{P})\ge\dd_n/2$.
Without loss of generality we may assume that $\hbox{dist}_4(\vec{u}_n, \mathcal{P})=\hbox{dist}_4(u_{n, 1}, \mathcal{P}_1)$.
Write $K(\vec{u}_n)=\vec{w}_n=(w_{n, 1}, w_{n, 2})$ and $w_{n, i}=\al_{n, i}\tilde{w}_{n, i}$ as in
 Lemma \ref{lemma4}. Then by the proof of Lemma \ref{lemma4}, we see that $\al_{n, i}$ are all uniformly bounded.
Combining this with (\ref{eq22-4}) and (\ref{eq2-13}), we deduce that
{\allowdisplaybreaks
\begin{align*}
\hbox{dist}_4(w_{n, 1}, \mathcal{P}_1)|w_{n, 1}^-|_4 &=|w_{n, 1}^-|_4^2\le C\|w_{n, 1}^-\|_{\la_1}^2\\
&\le C\int_{\Om}\left(|\nabla w_{n, 1}^-|^2+\la_1 (w_{n, 1}^-)^2-\bb t_2(\vec{u}_n)u_{n, 2}^2 (w_{n, 1}^-)^2\right)\,dx\\
&=-C\al_{n, 1}\mu_1 t_1(\vec{u}_n)\int_{\Om}u_{n, 1}^3 w_{n, 1}^-\,dx\\
&\le C\int_{\Om}(u_{n, 1}^-)^3 w_{n, 1}^-\,dx\le C|u_{n, 1}^-|_4^3|w_{n, 1}^-|_4\\
&=C\hbox{dist}_4(u_{n, 1}, \mathcal{P}_1)^3 |w_{n, 1}^-|_4
\le C \dd_n^3|w_{n, 1}^-|_4.
\end{align*}
}%
So $\hbox{dist}_4(K(\vec{u}_n), \mathcal{P})\le \hbox{dist}_4(w_{n, 1}, \mathcal{P}_1)\le C\dd_n^3<\dd_n/2$
holds for $n$ sufficiently large, which is a contradiction. This completes the proof.\hfill$\square$\\

Now let us define a map
$V: \mathcal{N}_b^\ast \to H$ by $ V(\vec{u}):=\vec{u}-K(\vec{u})$.
We will prove that $(\sqrt{t_1(\vec{u})}u_1, \sqrt{t_2(\vec{u})}u_2)$ is a sign-changing solution of (\ref{eq2})
if $\vec{u}=(u_1, u_2)\in\mathcal{N}_b\setminus\mathcal{P}$ satisfies $V(\vec{u})=0$.

\bl\label{lemma9} Let $\vec{u}_n=(u_{n, 1}, u_{n, 2})\in\mathcal{N}_b$ be such that
$$J_\bb(\vec{u}_n)\to c<d_k\quad\hbox{and}\quad V(\vec{u}_n)\to 0\quad\hbox{strongly in $H$}.$$
Then up to a subsequence, there exists $\vec{u}\in \mathcal{N}_b$ such that $\vec{u}_n\to \vec{u}$ strongly in $H$ and $V(\vec{u})=0$. \el

\noindent {\bf Proof. } By Lemma \ref{lemma4}, up to a subsequence, we may assume that $\vec{u}_n\rightharpoonup \vec{u}=(u_1, u_2)$ weakly in $H$ and
$\vec{w}_n:=K(\vec{u}_n)=(w_{n, 1}, w_{n, 2})\to \vec{w}=(w_1, w_2)$ strongly in $H$. Recalling $V(\vec{u}_n)\to 0$, we get
{\allowdisplaybreaks
\begin{align*}
&\int_{\Om}\nabla u_{n, i}\nabla (u_{n, i}-u_i)\,dx=\int_{\Om}\nabla (w_{n, i}-w_{i})\nabla (u_{n, i}-u_i)\,dx\\
&+\int_{\Om}\nabla w_{i}\nabla (u_{n, i}-u_i)\,dx
+\int_{\Om}\nabla (u_{n, i}-w_{n, i})\nabla (u_{n, i}-u_i)\,dx=o(1).
\end{align*}
}%
Then it is easy to see that $\vec{u}_n\to \vec{u}$ strongly in $H$ and so $\vec{u}\in \overline{\mathcal{N}_b}$. Hence $V(\vec{u})=\lim_{n\to\iy} V(\vec{u}_n)=0$. Moreover, $J_\bb(\vec{u})=c<d_k$ and so $\vec{u}\in\mathcal{N}_b$.\hfill$\square$

\bl\label{lemma10}Recall $C_0>0$ in Lemma \ref{lemma2}. Then
$$J_\bb'(\vec{u})[V(\vec{u})]\ge \frac{C_0}{2}\|V(\vec{u})\|_H^2,\quad \hbox{for any}\,\, \vec{u}\in\mathcal{N}_b.$$
\el

\noindent {\bf Proof. } Fix any $\vec{u}=(u_1, u_2)\in\mathcal{N}_b$ and write $\vec{w}=K(\vec{u})=(w_1, w_2)$ as above,
then $V(\vec{u})=(u_1-w_1, u_2-w_2)$. By (\ref{eq2-13}) we have
$\int_{\Om}u_i^3(u_i-w_i)\,dx=1-1=0$.
Then we deduce from (\ref{eq2-11})-(\ref{eq2-11-1}), (\ref{eq22-4}) and (\ref{eq2-13}) that
{\allowdisplaybreaks
\begin{align*}
&J_\bb'(\vec{u})[V(\vec{u})]\\
=&\sum_{i=1}^2t_i(\vec{u})\int_{\Om}\left(\nabla u_i\nabla (u_i-w_i)+\la_i u_i(u_i-w_i)-t_j(\vec{u})\bb u_i(u_i-w_i)u_j^2\right)\,dx\\
= &\sum_{i=1}^2t_i(\vec{u})\int_{\Om}\big(\nabla u_i\nabla (u_i-w_i)+\la_i u_i(u_i-w_i)-t_j(\vec{u})\bb w_i(u_i-w_i)u_j^2\\
&\qquad\qquad\qquad-t_j(\vec{u})\bb(u_i-w_i)^2u_j^2\big)\,dx\\
=&\sum_{i=1}^2t_i(\vec{u})\int_{\Om}\big(\nabla u_i\nabla (u_i-w_i)+\la_i u_i(u_i-w_i)-\nabla w_i\nabla(u_i-w_i)\\
&\qquad-\la_i w_i(u_i-w_i)+\al_i\mu_i t_i(\vec{u})u_i^3(u_i-w_i)-t_j(\vec{u})\bb(u_i-w_i)^2u_j^2\big)\,dx\\
=&\sum_{i=1}^2t_i(\vec{u})\int_{\Om}\big(|\nabla (u_i-w_i)|^2+\la_i |u_i-w_i|^2-t_j(\vec{u})\bb(u_i-w_i)^2u_j^2\big)\,dx\\
\ge&\sum_{i=1}^2\frac{t_i(\vec{u})}{2}\|u_i-w_i\|_{\la_i}^2\ge \frac{C_0}{2}\|V(\vec{u})\|_H^2.
\end{align*}
}%
This completes the proof.\hfill$\square$

\bl\label{lemma11} There exists a unique global solution $\eta=(\eta_1, \eta_2) : [0, \iy)\times \mathcal{N}_{b, \bb} \to H$ for the initial value problem
\be\label{eq2-19}\frac{d}{dt}\eta(t,\vec{u})=-V(\eta(t, \vec{u})), \quad \eta(0, \vec{u})=\vec{u}\in\mathcal{N}_{b,\bb}.\ee
Moreover,

\begin{itemize}
\item[$(i)$] $\eta(t, \vec{u})\in \mathcal{N}_{b, \bb}$ for any $ t>0$ and $\vec{u}\in\mathcal{N}_{b, \bb}$.
\item[$(ii)$] $\eta(t, \sg_i(\vec{u}))=\sg_i(\eta(t, \vec{u}))$ for any $t>0$, $\vec{u}\in\mathcal{N}_{b,\bb}$ and $i=1, 2$.
\item[$(iii)$] For every $\vec{u}\in\mathcal{N}_{b, \bb}$, the map $t\mapsto J_\bb(\eta(t, \vec{u}))$ is non-increasing.
\item[$(iv)$] There exists $\dd_0\in (0, 2^{-1/4})$ such that, for every $\dd<\dd_0$, there holds
$$\eta(t, \vec{u})\in \mathcal{P}_\dd \quad\hbox{whenever}\,\,\vec{u}\in \mathcal{N}_{b, \bb}\cap \mathcal{P}_\dd\,\,\hbox{and}\,\,t>0.$$
\end{itemize}

\el

\noindent {\bf Proof. } Recalling Lemma \ref{lemma3}, we have $V(\vec{u})\in C^1(\mathcal{N}_b^\ast, H)$.
Since $\mathcal{N}_{b, \bb}\subset\mathcal{N}_b^\ast$ and $\mathcal{N}_b^\ast$ is open, so (\ref{eq2-19}) has a unique solution
$\eta: [0, T_{\max})\times\mathcal{N}_{b, \bb}\to H$,
where $T_{\max}>0$ is the maximal time such that $\eta(t, \vec{u})\in \mathcal{N}_b^\ast$ for all
$t\in [0, T_{\max})$ (Note that $V(\cdot)$ is defined only on $\mathcal{N}_b^\ast$). We should prove $T_{\max}=+\iy$ for any $\vec{u}\in \mathcal{N}_{b, \bb}$.
Fixing any $\vec{u}=(u_1, u_2)\in \mathcal{N}_{b, \bb}$, we have
{\allowdisplaybreaks
\begin{align*}
\frac{d}{dt}\int_{\Om}\eta_i(t, \vec{u})^4\,dx &=-4\int_{\Om}\eta_i(t, \vec{u})^3(\eta_i(t, \vec{u})-K_i(\eta(t, \vec{u})))\,dx\\
&=4-4\int_{\Om}\eta_i(t, \vec{u})^4\,dx,\quad \forall\, 0< t<T_{\max}.
\end{align*}
}%
Recalling $\int_{\Om}\eta_i(0, \vec{u})^4\,dx=\int_{\Om}u_i^4\,dx=1$, we deduce that
$\int_{\Om}\eta_i(t, \vec{u})^4\,dx\equiv 1$ for all $0\le t<T_{\max}$.
So $\eta(t, \vec{u})\in \mathcal{M}$, that is $\eta(t, \vec{u})\in\mathcal{M}\cap\mathcal{N}_b^\ast=\mathcal{N}_b$ for all $t\in [0, T_{\max})$.
Assume by contradiction that $T_{\max}<+\iy$, then $\eta(T_{\max}, \vec{u})\in\partial\mathcal{N}_b$, and so
$J_\bb (\eta(T_{\max}, \vec{u}))\ge d_k$.
Since $\eta(t,\vec{u})\in \mathcal{N}_b$ for any $t\in [0, T_{\max})$, we deduce from Lemma \ref{lemma10} that
{\allowdisplaybreaks
\begin{equation}\label{eq2-20}
\begin{split}
J_\bb\left(\eta\left(T_{\max}, \vec{u}\right)\right)
&=J_\bb(\vec{u})-\int_0^{T_{\max}} J_\bb'(\eta(t, \vec{u}))[V(\eta(t, \vec{u}))]\,dt\\
&\le J_\bb(\vec{u})-\frac{C_0}{2}\int_0^{T_{\max}} \|V(\eta(t, \vec{u}))\|_H^2\,dt\le J_\bb(\vec{u})<d_k,
\end{split}
\end{equation}
}%
a contradiction. So $T_{\max}=+\iy$. Then similarly as (\ref{eq2-20}) we have $J_\bb(\eta(t, \vec{u}))\le J_\bb(\vec{u})<d_k$ for all $t>0$, so
$\eta(t, \vec{u})\in\mathcal{N}_{b, \bb}$ and then $(i), (iii)$ hold.

By (\ref{eq2-15}) we have $V(\sg_i(\vec{u}))=\sg_i(V(\vec{u}))$. Then by the uniqueness of solutions of the initial value
problem (\ref{eq2-19}), it is easy to check that
$(ii)$ holds.

Finally, let $\dd_0\in (0, 2^{-1/4})$ such that Lemma \ref{lemma8} holds for every $\dd<\dd_0$. For any $\vec{u}\in\mathcal{N}_{b, \bb}$ with $\hbox{dist}_4(\vec{u}, \mathcal{P})=\dd<\dd_0$, since
$$\eta(t, \vec{u})=\vec{u}+t\frac{d}{dt}{\eta}(0, \vec{u})+o(t)=\vec{u}-t V(\vec{u})+o(t)=(1-t)\vec{u}+tK(\vec{u})+o(t),$$
so we see from Lemma \ref{lemma8} that
{\allowdisplaybreaks
\begin{align*}
\hbox{dist}_4(\eta(t, \vec{u}), \mathcal{P})&=\hbox{dist}_4((1-t)\vec{u}+tK(\vec{u})+o(t), \mathcal{P})\\
&\le(1-t)\hbox{dist}_4(\vec{u}, \mathcal{P})+t\hbox{dist}_4(K(\vec{u}), \mathcal{P})+o(t)\\
&\le(1-t)\dd+t\dd/2+o(t)<\dd
\end{align*}
}%
for $t>0$ sufficiently small. Hence $(iv)$ holds.\hfill$\square$\\

\noindent {\bf Proof of Theorem \ref{th1}. }

{\bf Step 1.} Fix any $k_1, k_2\in [2, k+1]$ and take any $\dd\in (0, \dd_0)$.  We prove that (\ref{eq2}) has a sign-changing solution $(\tilde{u}_1, \tilde{u}_2)\in H$ such that
$E_\bb (\tilde{u}_1, \tilde{u}_2)=c_{\bb, \dd}^{k_1, k_2}$.

Write $c_{\bb, \dd}^{k_1, k_2}$ simply by $c$ in this step. Recall that $c<d_k$. We claim that
there exists a sequence $\{\vec{u}_n :n\ge 1\}\subset \mathcal{N}_{b, \bb}$ such that
\be\label{eq2-21}J_\bb(\vec{u}_n)\to c,\,\, V(\vec{u}_n)\to 0\,\,\,\hbox{as $n\to\iy$,}\,\,\,
\hbox{and}\,\, \hbox{dist}_4(\vec{u}_n, \mathcal{P})\ge \dd,\,\,\,\forall\,n\in\mathbb{N}.\ee

If (\ref{eq2-21}) does not hold, there exists small $\e\in (0, 1)$ such that
$$\|V(\vec{u})\|_H^2\ge \e,\quad\forall\,u\in\mathcal{N}_{b, \bb},\,\,|J_\bb(\vec{u})-c|\le 2\e,\,\,\hbox{dist}_4(\vec{u}, \mathcal{P})\ge \dd.$$
Recalling the definition of $c$ in (\ref{eq2-17}), we see that there exists $A\in \Gamma_\bb^{(k_1, k_2)}$ such that
$$\sup_{A\setminus \mathcal{P}_\dd}J_\bb<c+\e.$$
Since $\sup_{A} J_\bb<d_k$, so $A\subset \mathcal{N}_{b, \bb}$.
Then we can consider $B=\eta(4/C_0, A)$, where $\eta$ is in Lemma \ref{lemma11} and $C_0$ is in Lemma \ref{lemma2}.
Lemma \ref{lemma11}-$(i)$ yields $B\subset\mathcal{N}_{b, \bb}$. By Lemma \ref{lemma5}-$(ii)$ and Lemma \ref{lemma11}-$(ii)$ we have $B\in \Gamma^{(k_1, k_2)}$. Again by Lemma \ref{lemma11}-$(iii)$, we have
$\sup_{B} J_\bb\le\sup_A J_\bb<d_k$, that is $B\in \Gamma_\bb^{(k_1, k_2)}$ and so
$\sup_{B\setminus \mathcal{P}_\dd}J_\bb\ge c$.
Then by Lemma \ref{lemma6} we can take $\vec{u}\in A$ such that $\eta(4/C_0, \vec{u})\in B\setminus \mathcal{P}_\dd$ and
$$c-\e\le \sup_{B\setminus \mathcal{P}_\dd}J_\bb-\e< J_\bb (\eta(4/C_0, \vec{u})).$$
Since $\eta(t, \vec{u})\in \mathcal{N}_{b, \bb}$ for any $t\ge 0$,
Lemma \ref{lemma11}-$(iv)$ yields $\eta(t, \vec{u})\not\in\mathcal{P}_\dd$ for any $t\in [0, 4/C_0]$.
In particular, $\vec{u}\not\in\mathcal{P}_\dd$ and so $J_\bb(\vec{u})< c+\e$. Then
for any $t\in [0, 4/C_0]$, we have
$$c-\e< J_\bb (\eta(4/C_0, \vec{u}))\le J_\bb (\eta(t, \vec{u}))\le J_\bb (\vec{u})< c+\e,$$
which implies $\|V(\eta(t, \vec{u}))\|_H^2\ge \e$ and
\begin{align*}
\frac{d}{dt}J_\bb (\eta(t, \vec{u}))=-J_\bb'(\eta(t, \vec{u}))[V(\eta(t, \vec{u}))]\le-\frac{C_0}{2}\|V(\eta(t, \vec{u}))\|_H^2\le-\frac{C_0}{2} \e
\end{align*}
for every $t\in [0, 4/C_0]$. Hence,
$$c- \e< J_\bb (\eta(4/C_0, \vec{u}))\le J_\bb(\vec{u})-\int_{0}^{4/C_0}\frac{C_0}{2}\e\,dt<c+\e-2\e=c-\e,$$
a contradiction. Therefore (\ref{eq2-21}) holds. By Lemma \ref{lemma9}, up to a subsequence,
there exists $\vec{u}=(u_1, u_2)\in \mathcal{N}_{b, \bb}$ such that $\vec{u}_n\to \vec{u}$ strongly in $H$ and $V(\vec{u})=0$, $J_\bb(\vec{u})=c=c_{\bb, \dd}^{k_1, k_2}$. Since $\hbox{dist}_4(\vec{u}_n, \mathcal{P})\ge \dd$, so
$\hbox{dist}_4(\vec{u}, \mathcal{P})\ge \dd$, which implies that both $u_1$ and $u_2$ are sign-changing.
Since $V(\vec{u})=0$, so $\vec{u}=K(\vec{u})$. Combining this with (\ref{eq2-13})-(\ref{eq2-14}), we see that $\vec{u}$ satisfies
\be\label{eq2-22-1}\begin{cases}
-\Delta u_1+\la_1 u_1=\al_1\mu_1 t_1(\vec{u})u_1^3+\bb t_2(\vec{u})u_2^2 u_1,\\
-\Delta u_2+\la_2 u_2=\al_2\mu_2 t_2(\vec{u})u_2^3+\bb t_1(\vec{u})u_1^2 u_2.
\end{cases}\ee
Recall that $|u_i|_4=1$ and $t_i(\vec{u})$ satisfies (\ref{eq2-4}). Multiplying (\ref{eq2-22-1}) by $u_i$ and integrating over $\Om$, we obtain that $\al_1=\al_2=1$.
Again by (\ref{eq2-22-1}), we see that $(\tilde{u}_1, \tilde{u}_2):=(\sqrt{t_1(\vec{u})}u_1, \sqrt{t_2(\vec{u})}u_2)$ is a sign-changing solution of the original problem (\ref{eq2}). Moreover, (\ref{eq2-5}) and (\ref{eq2-9}) yield
$E_\bb (\tilde{u}_1, \tilde{u}_2)=J_\bb (u_1, u_2)=c_{\bb, \dd}^{k_1, k_2}$.

{\bf Step 2.} We prove that (\ref{eq2}) has at least $k$ sign-changing solutions.

Assume by contradiction that (\ref{eq2}) has at most $k-1$ sign-changing solutions.
Fix any $k_2\in [2, k+1]$ and $\dd\in (0, \dd_0)$.
Since $\Gamma_{\bb}^{(k_1+1, k_2)}\subset\Gamma_{\bb}^{(k_1, k_2)}$, we have
\begin{align}\label{eq6-2}
c_{\bb,\dd}^{2, k_2}\le c_{\bb, \dd}^{3, k_2}\le\cdots\le c_{\bb, \dd}^{k, k_2}\le c_{\bb, \dd}^{k+1, k_2}<d_k.
\end{align}
Since $c_{\bb,\dd}^{k_1, k_2}$ is a sign-changing critical value of $E_\bb$ for each $k_1\in [2, k+1]$
(that is,
$E_\bb$ has a sign-changing critical point $\vec{u}$ with $E_\bb (\vec{u})=c_{\bb,\dd}^{k_1, k_2}$),
by (\ref{eq6-2}) and our assumption that
(\ref{eq2}) has at most $k-1$ sign-changing solutions, there exists some $2\le N_1\le k$ such that
\be\label{eq2-27}c_{\bb,\dd}^{N_1, k_2}= c_{\bb,\dd}^{N_1+1, k_2}=:\bar{c}<d_k.\ee
Define
\be\label{eq6-10}\mathcal{K}:=\{\vec{u}\in \mathcal{N}_b \,\,\,:\,\,\,\vec{u}\,\,\hbox{sign-changing},\,\,\, J_\bb(\vec{u})=\bar{c},\,\,\,V(\vec{u})=0 \}.\ee
Then $\mathcal{K}$ is finite. By (\ref{eq2-15}) one
has that $\sg_i(\vec{u})\in \mathcal{K}$ if $\vec{u}\in\mathcal{K}$, that is, $\mathcal{K}\subset\mathcal{F}$.
Hence there exist $k_0\le k-1$ and $\{\vec{u}_m : 1\le m\le k_0\}\subset \mathcal{K}$
such that
$$\mathcal{K}=\{\vec{u}_m,\,\sg_1(\vec{u}_m),\, \sg_2(\vec{u}_m),\,-\vec{u}_m\,\,:\,\,\,1\le m\le k_0\}.$$
Then there exist open neighborhoods $O_{\vec{u}_m}$ of $\vec{u}_m$ in $H$, such that
any two of $\overline{O_{\vec{u}_{m}}}, \,\sg_1(\overline{O_{\vec{u}_{m}}}), \,\sg_2(\overline{O_{\vec{u}_{m}}})$ and $-\overline{O_{\vec{u}_{m}}}$,
where $1\le m\le k_0$, are disjointed and
$$
\mathcal{K}\subset O:=\bigcup_{m=1}^{k_0}
O_{\vec{u}_m}\cup\sg_1(O_{\vec{u}_m})\cup\sg_2(O_{\vec{u}_m})\cup-O_{\vec{u}_m}.
$$
Define a continuous map $\tilde{f}: \overline{O}\to \R\setminus\{0\}$ by
$$\tilde{f}(\vec{u}):=\begin{cases}1,&\hbox{if}\,\, \vec{u}\in \bigcup_{m=1}^{k_0}\overline{O_{\vec{u}_m}}\cup\sg_2(\overline{O_{\vec{u}_m}}),\\
-1, &\hbox{if}\,\, \vec{u}\in\bigcup_{m=1}^{k_0} \sg_1(\overline{O_{\vec{u}_m}})\cup-\overline{O_{\vec{u}_m}}.\end{cases}$$
Then
$\tilde{f}(\sg_1(\vec{u}))=-\tilde{f}(\vec{u})$ and $\tilde{f}(\sg_2(\vec{u}))=\tilde{f}(\vec{u}).$
By Tietze's extension theorem, there exists $f\in C(H, \R)$ such that $f|_{O}\equiv\tilde{f}$.
Define
$$F(\vec{u}):=\frac{f(\vec{u})+f(\sg_2(\vec{u}))-f(\sg_1(\vec{u}))-f(-\vec{u})}{4},$$
then $F|_O\equiv\tilde{f}$,
$F(\sg_1(\vec{u}))=-F(\vec{u})$ and $F(\sg_2(\vec{u}))=F(\vec{u})$.
Define
$$\mathcal{K}_\tau:=\left\{\vec{u}\in \mathcal{N}_b : \inf_{\vec{v}\in \mathcal{K}}\|\vec{u}-\vec{v}\|_{H}<\tau\right\}.$$
Then we can take small $\tau>0$ such that $\mathcal{K}_{2\tau}\subset O$. Recalling $V(\vec{u})=0$ in $\mathcal{K}$
and $\mathcal{K}$ finite, we see that there exists $\widetilde{C}>0$ such that
\be\label{eq2-28}\|V(\vec{u})\|_{H}\le\widetilde{C},\quad\forall\,\,\vec{u}\in \overline{\mathcal{K}_{2\tau}}.\ee
For any $\vec{u}\in \mathcal{K}_{2\tau}$,
we have $F(\vec{u})=\tilde{f}(\vec{u})\neq 0$. That is
$F(\mathcal{K}_{2\tau})\subset\R\setminus\{0\}$.
By (\ref{eq6-10}) and Lemma \ref{lemma9} there exists small $\e\in (0, (d_k-\bar{c})/2)$ such that
\be\label{eq2-30}\|V(\vec{u})\|_H^2\ge \e,\quad\forall\,u\in\mathcal{N}_b\setminus(\mathcal{K}_{\tau}\cup \mathcal{P}_\dd)
\,\,\,\hbox{satisfying}\,\,\,|J_\bb(\vec{u})-\bar{c}|\le 2\e.\ee
Recalling $C_0$ in Lemma \ref{lemma2}, we let
\be\label{eq2-32}\al:=\frac{1}{2}\min\left\{1, \frac{\tau C_0}{2\widetilde{C}}\right\}.\ee
By (\ref{eq2-17})-(\ref{eq2-17-1}) and (\ref{eq2-27})
we take $A\in \Gamma_{\bb}^{(N_1+1, k_2)}$ such that
\be\label{eq2-33}\sup_{A\setminus \mathcal{P}_{\dd}}J_\bb<c_{\bb,\dd}^{N_1+1, k_2}+\al\e/2=\bar{c}+\al\e/2.\ee
Let $B:=A\setminus \mathcal{K}_{2\tau}$, then it is easy to check that $B\subset\mathcal{F}$. We claim that $\vec{\ga}(B)\ge (N_1, k_2)$. If not, there exists
$\tilde{g}\in F_{(N_1, k_2)}(B)$ such that $\tilde{g}(\vec{u})\neq 0$ for any $\vec{u}\in B$.
By Tietze's extension theorem, there exists $\bar{g}=(\bar{g}_1, \bar{g}_2)\in C(H, \R^{N_1-1}\times\R^{k_2-1})$ such that $\bar{g}|_{B}\equiv\tilde{g}$.
Define $g=(g_1, g_2)\in C(H, \R^{N_1-1}\times\R^{k_2-1})$ by
{\allowdisplaybreaks
\begin{align*}g_1(\vec{u}):=\frac{\bar{g}_1(\vec{u})+\bar{g}_1(\sg_2(\vec{u}))-\bar{g}_1(\sg_1(\vec{u}))-\bar{g}_1(-\vec{u})}{4},\\
g_2(\vec{u}):=\frac{\bar{g}_2(\vec{u})+\bar{g}_2(\sg_1(\vec{u}))-\bar{g}_2(\sg_2(\vec{u}))-\bar{g}_2(-\vec{u})}{4},
\end{align*}
}%
then $g|_{B}\equiv\tilde{g}$,
$g_i(\sg_i(\vec{u}))=-g_i(\vec{u})$ and $g_i(\sg_j(\vec{u}))=g_i(\vec{u})$ for $j\neq i$.
Finally we define $G=(G_1, G_2)\in C(A,\, \R^{N_1+1-1}\times\R^{k_2-1})$ by
$$G_1(\vec{u}):=(F(\vec{u}),\, g_1(\vec{u}))\in \R^{N_1+1-1},\quad G_2(\vec{u}):=g_2(\vec{u})\in\R^{k_2-1}.$$
By our constructions of $F$ and $g$, we have $G\in F_{(N_1+1, k_2)}(A)$. Since $\vec{\ga}(A)\ge (N_1+1, k_2)$, so $G(\vec{u})=0$ for some $\vec{u}\in A$.
If $\vec{u}\in \mathcal{K}_{2\tau}$, then $F(\vec{u})\neq 0$, a contradiction. So $\vec{u}\in A\setminus\mathcal{K}_{2\tau}=B$,
and then $g(\vec{u})=\tilde{g}(\vec{u})\neq 0$, also a contradiction. Hence $\vec{\ga}(B)\ge (N_1, k_2)$.
Note that $B\subset A\subset \mathcal{N}_b$ and $\sup_{B}J_\bb\le \sup_{A}J_\bb<d_k$,
we see that $B\subset\mathcal{N}_{b, \bb}$ and $B\in \Gamma_{\bb}^{(N_1, k_2)}$.
Then we can consider $D:=\eta(\tau/(2\widetilde{C}), B)$, where $\eta$ is in Lemma \ref{lemma11} and $\widetilde{C}$ is in (\ref{eq2-28}). By Lemma \ref{lemma5}-$(ii)$ and Lemma \ref{lemma11} we have $D\subset \mathcal{N}_{b, \bb}$, $D\in \Gamma^{(N_1, k_2)}$ and
$\sup_{D} J_\bb\le\sup_B J_\bb<d_k$, that is $D\in \Gamma_{\bb}^{(N_1, k_2)}$.
Then we see from (\ref{eq2-17})-(\ref{eq2-17-1}) and (\ref{eq2-27}) that
$$\sup_{D\setminus \mathcal{P}_{\dd}}J_\bb\ge c_{\bb,\dd}^{N_1, k_2}=\bar{c}.$$
By Lemma \ref{lemma6} we can take $\vec{u}\in B$ such that $\eta(\tau/(2\widetilde{C}), \vec{u})\in D\setminus \mathcal{P}_{\dd}$ and
$$\bar{c}-\al\e/2\le \sup_{D\setminus \mathcal{P}_\dd}J_\bb-\al\e/2< J_\bb (\eta(\tau/(2\widetilde{C}), \vec{u})).$$
Since $\eta(t, \vec{u})\in \mathcal{N}_{b, \bb}$ for any $t\ge 0$, Lemma \ref{lemma11}-$(iv)$ yields $\eta(t, \vec{u})\not\in\mathcal{P}_\dd$ for any $t\in [0, \tau/(2\widetilde{C})]$. In particular, $\vec{u}\not\in\mathcal{P}_{\dd}$ and so (\ref{eq2-33}) yields $J_\bb(\vec{u})< \bar{c}+\al\e/2$. Then
for any $t\in [0, \tau/(2\widetilde{C})]$, we have
$$\bar{c}-\al\e/2< J_\bb (\eta(\tau/(2\widetilde{C}), \vec{u}))\le J_\bb (\eta(t, \vec{u}))\le J_\bb (\vec{u})< \bar{c}+\al\e/2.$$
Recall that $\vec{u}\in B=A\setminus\mathcal{K}_{2\tau}$. If there exists $T\in (0, \tau/(2\widetilde{C}))$ such
that $\eta(T, \vec{u})\in \mathcal{K}_\tau$, then
there exist $0\le t_1<t_2\le T$ such that $\eta(t_1, \vec{u})\in\partial\mathcal{K}_{2\tau}$, $\eta(t_2, \vec{u})\in\partial\mathcal{K}_{\tau}$ and
$\eta(t, \vec{u})\in\mathcal{K}_{2\tau}\setminus\mathcal{K}_\tau$ for any $t\in (t_1, t_2)$. So we see from (\ref{eq2-28}) that
$$\tau\le\|\eta(t_1, \vec{u})-\eta(t_2, \vec{u})\|_H=\left\|\int_{t_1}^{t_2}V(\eta(t, \vec{u}))\,dt\right\|_H\le 2\widetilde{C}(t_2-t_1),$$
that is, $\tau/(2\widetilde{C})\le t_2-t_1\le T$, a contradiction. Hence $\eta(t, \vec{u})\not\in \mathcal{K}_\tau$ for any $t\in (0, \tau/(2\widetilde{C}))$.
Then as Step 1, we deduce from (\ref{eq2-30}) and (\ref{eq2-32}) that
\begin{align*}
\bar{c}-\frac{\al\e}{2}< J_\bb (\eta(\tau/(2\widetilde{C}), \vec{u}))\le J_\bb(\vec{u})-\int_0^{\frac{\tau}{2\widetilde{C}}}\frac{C_0}{2}\e\,dt<\bar{c}+\frac{\al\e}{2}-\al\e=\bar{c}-\frac{\al\e}{2},
\end{align*}
a contradiction.
Hence (\ref{eq2}) has at least $k$ sign-changing solutions for any $\bb\in (0, \bb_k)$.
This completes the proof.\hfill$\square$

\vskip0.1in

\s{Proof of Theorem \ref{th2}}
\renewcommand{\theequation}{3.\arabic{equation}}

Let $k=1$ in Section 2. By the proof of Theorem \ref{th1} there exists $\bb_1>0$ such that, for any $\bb\in (0, \bb_1)$, (\ref{eq2}) has a sign-changing solution $(u_{\bb, 1}, v_{\bb, 1})$ with $E_\bb(u_{\bb, 1}, v_{\bb, 1})=c_{\bb, \dd}^{2, 2}<d_1$. Recalling $S$ in (\ref{eq2-10}), we define
\be\label{eq3-1}\bb_1':=\min\left\{S^2/(4d_1),\,\bb_1\right\}.\ee
Fix any $\bb\in (0, \bb_1')$ and define
$$c_\bb:=\inf_{\vec{u}\in\mathcal{K}_\bb}E_\bb(\vec{u});\quad \mathcal{K}_\bb:=\{\vec{u} : \vec{u}\,\,\hbox{is a sign-changing solution of (\ref{eq2})}\}.$$
Then $\mathcal{K}_\bb\neq \emptyset$ and $c_\bb<d_1$. Let $\vec{u}_n=(u_{n, 1}, u_{n, 2})\in\mathcal{K}_\bb$ be a minimizing sequence of $c_\bb$ with
$E_\bb(\vec{u}_n)<d_1$ for all $n\ge 1$. Then $\|u_{n, 1}\|_{\la_1}^2+\|u_{n, 2}\|_{\la_2}^2<4d_1$. Up to a subsequence, we may assume that $\vec{u}_n\to \vec{u}=(u_1, u_2)$ weakly in $H$ and strongly in $L^4(\Om)\times L^4(\Om)$. Since $E_\bb'(\vec{u}_n)=0$, it is standard to prove that $\vec{u}_n\to \vec{u}=(u_1, u_2)$ strongly in $H$, $E_\bb'(\vec{u})=0$ and $E_\bb(\vec{u})=c_\bb$. On the other hand, we deduce
from $E_\bb'(\vec{u}_n)(u_{n, 1}^\pm, 0)=0$ and $E_\bb'(\vec{u}_n)(0, u_{n, 2}^\pm)=0$ that
\begin{align*}
S|u_{n, i}^{\pm}|_4^2&\le\|u_{n, i}^{\pm}\|_{\la_i}^2=\mu_i|u_{n, i}^{\pm}|_4^4+\bb\int_{\Om}|u_{n, i}^{\pm}|^2 u_{n, j}^2dx
\le \mu_i|u_{n, i}^{\pm}|_4^4+\bb|u_{n, i}^{\pm}|_4^2 |u_{n, j}|_4^2\\
&\le\mu_i|u_{n, i}^{\pm}|_4^4+\frac{\bb}{S}|u_{n, i}^{\pm}|_4^2\|u_{n, j}\|_{\la_j}^2<\mu_i|u_{n, i}^{\pm}|_4^4+\frac{4d_1\bb}{S}|u_{n, i}^{\pm}|_4^2,
\end{align*}
which implies that $|u_{n, i}^{\pm}|_4\ge C>0$ for all $n\ge 1$ and $i=1, 2$, where $C$ is a constant independent of $n$ and $i$.
Hence $|u_{i}^\pm|_4\ge C$ and so $\vec{u}$ is a least energy sign-changing solution of (\ref{eq2}).\hfill$\square$

\vskip0.1in

\s{Proof of Theorems \ref{th3}}
\renewcommand{\theequation}{4.\arabic{equation}}

The following arguments are similar to those in Section 2 with some modifications. Here, although some definitions
are slight different from those in Section 2, we will use the same notations as in Section 2 for convenience.
To obtain semi-nodal solutions $(u_1, u_2)$ such that $u_1$ changes sign and $u_2$ is positive,
we consider the following functional
\begin{align*}
\widetilde{E}_\bb(u_1, u_2):=\frac{1}{2} \left(\|u_1\|_{\la_1}^2+\|u_2\|_{\la_2}^2\right)
-\frac{1}{4}\left(\mu_1 |u_1|_4^4+\mu_2|u_2^+|_4^4\right)-\frac{\beta}{2} \int_{\Om}u_1^2u_2^2\,dx,
\end{align*}
and modify the definition of $\wH$ by
$\wH:=\{(u_1, u_2)\in H : u_1\neq 0,\,\,u_2^+\neq 0\}$.
Then by similar proofs as in Section 2, we have the following lemma.

\bl\label{lemma02}For any $\vec{u}=(u_1, u_2)\in\wH$, if
\be\label{eq02-3}\begin{cases}
\mu_2|u_2^+|_4^4\|u_1\|_{\la_1}^2-\bb\|u_2\|_{\la_2}^2\int_{\Om}u_1^2u_2^2\,dx>0,\\
\mu_1|u_1|_4^4\|u_2\|_{\la_2}^2-\bb\|u_1\|_{\la_1}^2\int_{\Om}u_1^2u_2^2\,dx>0,
\end{cases}\ee
then system
\be\label{eq02-4}
\begin{cases}\|u_1\|_{\la_1}^2 =
t_1\mu_1 |u_1|_4^4+t_2\beta\int_{\Om} u_1^2 u_2^2\,dx\\
\|u_2\|_{\la_2}^2 =
t_2\mu_2 |u_2^+|_4^4+t_1\beta\int_{\Om} u_1^2 u_2^2\,dx\end{cases}\ee
has a unique solution
\be\label{eq02-5}
\begin{cases}t_1(\vec{u})=\frac{\mu_2|u_2^+|_4^4\|u_1\|_{\la_1}^2-\bb\|u_2\|_{\la_2}^2\int_{\Om}u_1^2u_2^2\,dx}{\mu_1\mu_2|u_1|_4^4 |u_2^+|_4^4-\bb^2(\int_{\Om}u_1^2 u_2^2\,dx)^2}>0\\
t_2(\vec{u})=\frac{\mu_1|u_1|_4^4\|u_2\|_{\la_2}^2-\bb\|u_1\|_{\la_1}^2\int_{\Om}u_1^2u_2^2\,dx}{\mu_1\mu_2|u_1|_4^4 |u_2^+|_4^4-\bb^2(\int_{\Om}u_1^2 u_2^2)^2\,dx}>0.\end{cases}\ee
Moreover,
{\allowdisplaybreaks
\begin{align}\label{eq02-6}\sup_{t_1, t_2\ge 0}&\widetilde{E}_\bb \left(\sqrt{t_1}u_1, \,\sqrt{t_2} u_2\right)=\widetilde{E}_\bb\left(\sqrt{t_1(\vec{u})}u_1, \sqrt{t_2(\vec{u})}u_2\right)\nonumber\\
&=\frac{1}{4}\frac{\mu_2|u_2^+|_4^4\|u_1\|_{\la_1}^4
-2\bb\|u_1\|_{\la_1}^2\|u_2\|_{\la_2}^2\int_{\Om}u_1^2u_2^2\,dx+\mu_1|u_1|_4^4\|u_2\|_{\la_2}^4}{\mu_1\mu_2|u_1|_4^4 |u_2^+|_4^4-\bb^2(\int_{\Om}u_1^2 u_2^2\,dx)^2}
\end{align}
}%
and $(t_1(\vec{u}), t_2(\vec{u}))$ is the unique maximum point of $\widetilde{E}_\bb (\sqrt{t_1}u_1, \sqrt{t_2}u_2)$.
\el

Now we modify the definitions of $\mathcal{M}^\ast$, $\mathcal{M}_\bb^\ast$, $\mathcal{M}_\bb^{\ast\ast}$, $\mathcal{M}$ and $\mathcal{M}_\bb$ by
{\allowdisplaybreaks
\begin{align}
\label{eq02-8-1}&\mathcal{M}^\ast:=\left\{\vec{u}\in H\,\,:\,\, 1/2<|u_1|^4_4<2,\,\, 1/2<|u_2^+|^4_4<2\right\};\\
&\mathcal{M}_\bb^\ast:=\left\{\vec{u}\in \mathcal{M}^\ast\,\,:\,\, \hbox{$\vec{u}$ satisfies (\ref{eq02-3})}\right\};\nonumber\\
&\mathcal{M}_\bb^{\ast\ast}:=\left\{\vec{u}\in \mathcal{M}^\ast\,\,:\,\, \begin{array}{ll}
\mu_2\|u_1\|_{\la_1}^2-\bb\|u_2\|_{\la_2}^2\int_{\Om}u_1^2u_2^2\,dx>0\\
\mu_1\|u_2\|_{\la_2}^2-\bb\|u_1\|_{\la_1}^2\int_{\Om}u_1^2u_2^2\,dx>0
\end{array}\right\};\nonumber\\
\label{eq02-8}&\mathcal{M}:=\left\{\vec{u}\in H\,\,:\,\, |u_1|_4=1,\,\, |u_2^+|_4=1\right\},\quad \mathcal{M}_\bb:=\mathcal{M}\cap \mathcal{M}_\bb^\ast,
\end{align}
}%
and define a new functional $J_\bb : \mathcal{M}_\bb^{\ast\ast}\to (0, +\iy)$ as in Section 2 by
$$J_\bb (\vec{u}):=\frac{1}{4}\frac{\mu_2\|u_1\|_{\la_1}^4
-2\bb\|u_1\|_{\la_1}^2\|u_2\|_{\la_2}^2\int_{\Om}u_1^2u_2^2\,dx+\mu_1\|u_2\|_{\la_2}^4}{\mu_1\mu_2-\bb^2(\int_{\Om}u_1^2 u_2^2\,dx)^2}
.$$
Then
$J_\bb\in C^1(\mathcal{M}^{\ast\ast}_\bb, \,(0, +\iy))$
and (\ref{eq2-11})-(\ref{eq2-11-1})
hold for any $\vec{u}\in \mathcal{M}_\bb$ and $\vp,\, \psi\in H_0^1(\Om)$.
Moreover,
Lemma \ref{lemma02} yields
\be\label{eq02-9}J_\bb (u_1, u_2)=\sup_{t_1, t_2\ge 0}\widetilde{E}_\bb \left(\sqrt{t_1}u_1, \,\sqrt{t_2} u_2\right),\quad \forall\,(u_1, u_2)\in \mathcal{M}_\bb.\ee

Under this new definitions (\ref{eq02-8-1})-(\ref{eq02-8}), we define $\mathcal{N}_b^\ast$ and $\mathcal{N}_b$ as in (\ref{eq22-1})-(\ref{eq22-3}).
Since $|u_2|_4^2\le S^{-1}\|u_2\|_{\la_2}^2\le b/S$ for all $\vec{u}\in\mathcal{N}_b^\ast$, by trivial modifications
it is easy to check that Lemmas \ref{lemma2} and \ref{lemma22} also hold here. Moreover,
we may assume that (\ref{eq22-4}) also holds here for any $\bb\in (0, \bb_k)$.

Now we fix any $\bb\in (0, \bb_k)$. For any $\vec{u}=(u_1, u_2)\in \mathcal{N}_b^*$, let $\tilde{w}_i\in H_0^1(\Om)$, $i=1, 2$, be the unique
solutions of the following linear problem
\be\label{eq02-012}\begin{cases}
-\Delta \tilde{w}_1+\la_1 \tilde{w}_1-\bb t_2(\vec{u})u_2^2 \tilde{w}_1=\mu_1 t_1(\vec{u})u_1^3,\quad \tilde{w}_1\in H_0^1(\Om),\\
-\Delta \tilde{w}_2+\la_2 \tilde{w}_2-\bb t_1(\vec{u})u_1^2 \tilde{w}_2=\mu_2 t_2(\vec{u})(u_2^+)^3,\quad \tilde{w}_2\in H_0^1(\Om).
\end{cases}\ee
As in Section 2, we define
\be\label{eq02-12}w_i=\al_i \tilde{w}_i,\quad\hbox{where}\,\,\al_1=\frac{1}{\int_{\Om}u_1^3 \tilde{w}_1\,dx}>0,\,\,
\al_2=\frac{1}{\int_{\Om}(u_2^+)^3 \tilde{w}_2\,dx}>0 .\ee
Then $(w_1, w_2)$ is the unique solution of the problem
\be\label{eq02-13}\begin{cases}
-\Delta w_1+\la_1 w_1-\bb t_2(\vec{u})u_2^2 w_1=\al_1\mu_1 t_1(\vec{u})u_1^3,\quad w_1\in H_0^1(\Om),\\
-\Delta w_2+\la_2 w_2-\bb t_1(\vec{u})u_1^2 w_2=\al_2\mu_2 t_2(\vec{u})(u_2^+)^3,\quad w_2\in H_0^1(\Om),\\
\int_{\Om}u_1^3 w_1\,dx=1,\quad \int_{\Om}(u_2^+)^3 w_2\,dx=1.
\end{cases}\ee
As in Section 2, the operator $K=(K_1, K_2) : \mathcal{N}_b^\ast\to H$ is defined as
$K(\vec{u}):=\vec{w}=(w_1, w_2)$, and similar arguments as Lemma \ref{lemma3}
yield $K\in C^1(\mathcal{N}_b^*, H)$. Since
$u_n\to u$ in $L^4(\Om)$ implies $u_n^+\to u^+$ in $L^4(\Om)$,
so Lemma \ref{lemma4} also holds for this new $K$ defined here. Clearly
\be\label{eq02-15}K(\sg_1(\vec{u}))=\sg_1(K(\vec{u})).\ee
 Remark that (\ref{eq02-15}) only holds for $\sg_1$ and in the sequel we only use $\sg_1$.
Consider
$$\mathcal{F}=\{A\subset \mathcal{M} : A\,\,\hbox{ is closed and}\,\,\sg_1(\vec{u})\in A\,\,\,\forall\,\,\vec{u}\in A\},$$
and, for each $A\in \mathcal{F}$ and $k_1\ge 2$, the class of functions
$$F_{(k_1, 1)}(A)=\left\{f: A\to \R^{k_1-1} \,:  \,f\,\,\hbox{continuous and}\,\,
 f(\sg_1(\vec{u}))=-f(\vec{u})
 \right\}.$$
\begin{definition} (Modified vector genus, slightly different from Definition \ref{definition1})
Let $A\in \mathcal{F}$ and take any $k_1\in\mathbb{N}$ with $k_1\ge 2$. We say that $\vec{\ga}(A)\ge (k_1, 1)$ if for every $f\in F_{(k_1, 1)}(A)$ there exists $\vec{u}\in A$ such that $f(\vec{u})=0$. We denote
$$\Gamma^{(k_1, 1)}:=\{A\in\mathcal{F} : \vec{\ga}(A)\ge (k_1, 1)\}.$$
\end{definition}

\bl\label{lemma05} (see \cite[Lemma 4.2]{CLZ}) With the previous notations, the following properties hold.
\begin{itemize}
\item[$(i)$] Take $A:=A_1\times A_2\subset \mathcal{M}$ and let $\eta: S^{k_1-1}\to A_1$ be a homeomorphism such that
$\eta(-x)=-\eta(x)$ for every $x\in S^{k_1-1}$. Then $A\in \Gamma^{(k_1, 1)}$.

\item[$(ii)$] We have $\overline{\eta(A)}\in \Gamma^{(k_1, 1)}$ whenever $A\in \Gamma^{(k_1, 1)}$ and a continuous map $\eta : A\to \mathcal{M}$ is such
that $\eta\circ \sg_1=\sg_1\circ\eta$.
\end{itemize}\el

Now we modify the definitions of
$\mathcal{P}$ and $\hbox{dist}_4(\vec{u}, \mathcal{P})$ in (\ref{cone})-(\ref{cone1}) by
\be\label{eq4-1}\mathcal{P}:=\mathcal{P}_1\cup -\mathcal{P}_1,\quad\hbox{dist}_4(\vec{u}, \mathcal{P}):=\min\big\{\hbox{dist}_4(u_1,\,\mathcal{P}_1),\,\,\hbox{dist}_4(u_1,\,-\mathcal{P}_1)\big\}.\ee
Under this new definition, $u_1$ changes sign if $\hbox{dist}_4(\vec{u}, \mathcal{P})>0$.

\bl\label{lemma16} (see \cite[Lemma 4.3]{CLZ}) Let $k_1\ge 2$. Then for any $\dd<2^{-1/4}$ and any $A\in\Gamma^{(k_1, 1)}$ there holds $A\setminus \mathcal{P}_\dd\neq\emptyset$.\el

\bl\label{lemma07} There exists $A\in\Gamma^{(k+1, 1)}$
such that $A\subset \mathcal{N}_b$ and $\sup_{A}J_\bb<d_k$.\el

\noindent {\bf Proof. } Recalling $\vp_0\in W_{k+1}$ is positive, we define
$$A_1:=\big\{u\in W_{k+1} \,:\, |u|_4=1\big\}, \quad  A_2:=\{C\vp_0 : C=1/|\vp_0|_4\}.$$
Then by Lemma \ref{lemma05}-$(i)$ one has $A:=A_1\times A_2\in \Gamma^{(k+1, 1)}$. The rest of the proof is the same as Lemma \ref{lemma7}.
\hfill$\square$\\

For every $k_1\in [2, k+1]$ and $0<\dd < 2^{-1/4}$, we define
$$ c_{\bb,\dd}^{k_1, 1}:=\inf_{A\in \Gamma_\bb^{(k_1, 1)}}\sup_{\vec{u}\in A\setminus \mathcal{P}_\dd}J_\bb(\vec{u}),$$
where the definition of $\Gamma_\bb^{(k_1, 1)}$ is the same as (\ref{eq2-17-1}).
Then Lemma \ref{lemma07} yields $\Gamma_\bb^{(k_1, 1)}\neq \emptyset$ and so $c_{\bb, \dd}^{k_1, 1}$ is well defined for each $k_1\in [2, k+1]$. Moreover,
$c_{\bb, \dd}^{k_1, 1}<d_k$ for any $\dd\in (0, 2^{-1/4})$ and $k_1\in [2, k+1]$. Define $\mathcal{N}_{b,\bb}:=\{\vec{u}\in  \mathcal{N}_b: J_\bb(\vec{u})<d_k\}$
as in Section 2.
Under the new definition (\ref{eq4-1}), it is easy to see that Lemma \ref{lemma8} also holds here.
Now as in Section 2, we define a map
$V : \mathcal{N}_b^\ast \to H$ by $V(\vec{u}):=\vec{u}-K(\vec{u}).$
Then Lemma \ref{lemma9} also holds here. Recall from (\ref{eq02-8}) and (\ref{eq02-13}) that
$\int_{\Om}(u_2^+)^3 (u_2-w_2)\,dx=1-1=0$ for any $\vec{u}=(u_1, u_2)\in \mathcal{N}_b$.
Then by similar arguments, we see that Lemma \ref{lemma10} also holds here.

\bl\label{lemma011} There exists a unique global solution $\eta=(\eta_1, \eta_2) : [0, \iy)\times \mathcal{N}_{b,\bb} \to H$ for the initial value problem
\be\label{eq02-19}\frac{d}{dt}\eta(t,\vec{u})=-V(\eta(t, \vec{u})), \quad \eta(0, \vec{u})=\vec{u}\in\mathcal{N}_{b, \bb}.\ee
Moreover, conclusions $(i)$, $(iii)$ and $(iv)$ of Lemma \ref{lemma11} also hold here, and
$\eta(t, \sg_1(\vec{u}))=\sg_1(\eta(t, \vec{u}))$ for any $t>0$ and $u\in\mathcal{N}_{b, \bb}$.

\el

\noindent {\bf Proof. }Recalling $V(\vec{u})\in C^1(\mathcal{N}_b^\ast, H)$, we see that
 (\ref{eq02-19}) has a unique solution
$\eta: [0, T_{\max})\times\mathcal{N}_{b, \bb}\to H$,
where $T_{\max}>0$ is the maximal time such that $\eta(t, \vec{u})\in \mathcal{N}_b^\ast$ for all
$t\in [0, T_{\max})$.
Fix any $\vec{u}=(u_1, u_2)\in \mathcal{N}_{b, \bb}$, we deduce from (\ref{eq02-19}) that
$
\frac{d}{dt}\int_{\Om}\left(\eta_2(t, \vec{u})^+\right)^4\,dx=4-4\int_{\Om}\left(\eta_2(t, \vec{u})^+\right)^4\,dx,\,\, \forall\, 0< t<T_{\max}.
$
Since $\int_{\Om}\left(\eta_2(0, \vec{u})^+\right)^4dx=\int_{\Om}(u_2^+)^4dx=1$, so
$\int_{\Om}\left(\eta_2(t, \vec{u})^+\right)^4 dx\equiv 1$
for all $0\le t<T_{\max}$.
Recalling (\ref{eq02-15}), we see that the rest of the proof is similar to Lemma \ref{lemma11}.\hfill$\square$\\

\noindent {\bf Proof of Theorem \ref{th3}. }
First we fix any $k_1\in [2, k+1]$.
Then by similar arguments as Step 1 in the proof of Theorem \ref{th1}, for small $\dd>0$,
there exists $\vec{u}=(u_1, u_2)\in\mathcal{N}_b$ such that
$$J_\bb(\vec{u})=c_{\bb,\dd}^{k_1, 1},\quad V(\vec{u})=0\quad\hbox{and dist}_4(\vec{u}, \mathcal{P})\ge \dd.$$
So $u_1$ changes sign. Since $V(\vec{u})=0$, so $\vec{u}=K(\vec{u})$. Combining this with (\ref{eq02-13}), we see that $\vec{u}$ satisfies
\be\label{eq02-22-1}\begin{cases}
-\Delta u_1+\la_1 u_1=\al_1\mu_1 t_1(\vec{u})u_1^3+\bb t_2(\vec{u})u_2^2 u_1,\\
-\Delta u_2+\la_2 u_2=\al_2\mu_2 t_2(\vec{u})(u_2^+)^3+\bb t_1(\vec{u})u_1^2 u_2.
\end{cases}\ee
Since $|u_1|_4=1$, $|u_2^+|_4=1$ and $t_i(\vec{u})$ satisfies (\ref{eq02-4}), so $\al_1=\al_2=1$.
Multiplying the second equation of (\ref{eq02-22-1}) by $u_2^-$ and integrating over $\Om$,
we see from (\ref{eq22-4}) that $\|u_2^-\|_{\la_2}^2=0$, so $u_2\ge 0$. By the strong
maximum principle, $u_2>0$ in $\Om$. Hence $(\tilde{u}_1, \tilde{u}_2):=(\sqrt{t_1(\vec{u})}u_1, \sqrt{t_2(\vec{u})}u_2)$ is a semi-nodal solution of the original problem (\ref{eq2}) with $\tilde{u}_1$ sign-changing and $\tilde{u}_2$ positive. Moreover, (\ref{eq02-6}) and (\ref{eq02-9}) yield
$E_\bb (\tilde{u}_1, \tilde{u}_2)=\widetilde{E}_\bb (\tilde{u}_1, \tilde{u}_2)=J_\bb (u_1, u_2)=c_{\bb, \dd}^{k_1, 1}< d_k$.
Finally, since $k_1\in [2, k+1]$,
by similar arguments as Step 2 of proving Theorem \ref{th1} with trivial modifications, we can prove that (\ref{eq2})
has at least $k$ semi-nodal solutions. This completes the proof.
\hfill$\square$

\br By a similar argument as in Section 3, we can prove that there exists $\bb_1''>0$ such that for any $\bb\in (0, \bb_1'')$, (\ref{eq2}) has a semi-nodal solution which has the least energy among all semi-nodal solutions.\er

\end{document}